\input epsf
\documentclass[12pt]{amsart}
\usepackage{amsfonts,amssymb,amscd}
\newcommand{\R}{\mathbb R}
\newcommand{\C}{\mathbb C}

\newcommand{\Z}{\mathbb Z}

\newtheorem{Theorem}{Theorem}[section]

\newtheorem{Corollary}{Corollary}[section]
\newtheorem{Proposition}{Proposition}[section]
\newtheorem{Remark}{Remark}[section]
\allowdisplaybreaks[2]
\numberwithin{equation}{section}
\numberwithin{figure}{section}
\title{Braids, their properties and generalizations}
\author[Vershinin]{V.~V.~Vershinin}
\address{D\'epartement des Sciences Math\'ematiques,
                                     Universit\'e Montpellier II,
Place Eug\'ene Bataillon,
34095 Montpellier cedex 5, France}
\email{ vershini@math.univ-montp2.fr}
\address{ Sobolev Institute of Mathematics, Novosibirsk, 630090,
Russia }
\email{ versh@math.nsc.ru}
\thanks{
The author was supported
 in part by the 
 by CNRS-NSF grant No~17149, INTAS grant No~03-5-3251
and the ACI project ACI-NIM-2004-243 "Braids and Knots".}
\begin{document}
\begin{abstract}
In the paper we give a survey on braid groups and subjects connected with them.
We start with the initial definition, then we give several interpretations
as well as several presentations of these groups. Burau presentation for the 
pure braid group and the Markov normal form are given next. Garside normal form
and his solution of the conjugacy problem are presented as well as more recent
results on the ordering and on the linearity of braid groups. Next topics are the 
generalizations
of braids, their homological properties and connections with the other mathematical
fields, like knot theory (via Alexander and Markov theorems) and homotopy groups 
of spheres.

\end{abstract}
\maketitle
\tableofcontents

\section{Introduction}

Braid groups describe intuitive concept of classes of continuous
deformations of braids, which are collections of interwining strands
whose endpoints are fixed.
Mathematically they can be 
 considered from various points of view. The first intuitive
approach is formalized  naturally as isotopy classes of a 
collection of $n$
connected curves (strings) in 3-dimensional space. 
This point
of view is connected with the definition of braid group as the
fundamental group of configuration space of $n$ points on a plane. 
Also braids can be interpreted as a mapping class group of 
a punctured disc and as a subgroup of automorphism group of a free
group (subsection \ref{subsec:aut}).

The present survey is organized as follows. 
  In Section~2 we make some historical remark.
Definition and general properties are considered in Section~3.
 Configuration spaces appear in Subsection~3.2.
Connections with groups of
automorphisms of free groups are given in Subsection~3.3
Presentations of the braid group which appeared quite recently
are observed in Subsection~3.4.
Section~\ref{sec:gar} is devoted to F.~A.~Garside's classical work
\cite{Gar} and Section~\ref{sec:dehor} to that of P.~Dehornoy
on ordering for braids. Representations and in particular linearity
are discussed in 
Section~\ref{sec:rep}.
In Section~\ref{sec:gen} various generalization of braids are
presented. Homological properties are observed in 
Section~\ref{sec:hom}. In the last Section~\ref{sec:con} we discuss 
connection with 
the knot theory given by the Alexander and Markov Theorems and 
with the homotopy groups of spheres.

The author is grateful to
F.~R.~Cohen and Wu Jie for fruitful discussions on theory of 
braids and their applications.
The author is also  thankful to Vik.~S.~Kulikov who attracted his attention
to the works of O.~Zariski \cite{Za1}, \cite{Za2} and
to E.~P.~Volokitin for the help and advices on the
presentation of the paper.

\section{Historical Remarks}
Braids were
rigorously defined by E.~Artin \cite{Art1} in 1925, although the roots
of this natural concept are seen in the works of A.~Hurwitz (\cite{Hu},
1891), R.~Fricke and F.~Klein (\cite{FK}, 1897) and even in  the
notebooks of K.-F.~Gauss. E.~Artin \cite{Art1} gave the presentation of the braid
group (see formulas (\ref{eq:brelations}) in Section~3) which is common now.
Already in the book of Felix Klein \cite{Kl} published in 1926 there 
appeared a chapter about braids. Essential topics about braids were also
presented in the Reidemeister's Knotentheorie \cite{R} published in 
1932.

In 30-ies there appeared a series of papers of Werner Burau \cite{Bu1}, 
\cite{Bu2}, \cite{Bu3} where he in particular gave the presentation 
of the pure braid group (see subsection \ref{subsec:pure}) and introduced 
the Burau representation (subsection \ref{subsec:bur}).\ 
Wilhelm Magnus in his work \cite{Mag1} published in 1934 established 
relations between braid groups and the mapping class groups.
At the same time 
there appeared the work of A.~A.~Markov \cite{Mar1}, which together with
Alexander Theorem \cite{Al} builds a bijection between links and equivalence 
classes of braids. It became an essential ingredient
in study of links and knots (in the work of
V.~F.~R.~Jones \cite{J}, for example). 
In 1936-37 were published the works of O.~Zariski \cite{Za1}, \cite{Za2},
where he 
discovered connections between braid groups and the fundamental group of 
the complement of discriminant
of the general polynomial
$$f_n(t)=a_0t^n+a_1t^{n-1}+...+a_{n-1}t+a_n ,$$
a point of view later rediscovered by V.~I.~Arnold \cite{Arn2}.
Zariski also understood connections between braids and configuration
spaces,
 gave the presentation of the braid group of the sphere, and
studied the braid groups of Riemann surfaces.
Amazingly and unfortunately these works of Zariski were not noticed 
by the specialists on braids and are not mentioned even in books and papers
where the presentations of braids of surfaces are discussed. 

In the beginning of 60-ies R. Fox and L. Neuwirth \cite{FoN} and 
E.~Fadell and L.~Neuwirth \cite{FaN} 
studied configuration spaces which turned out to be $K(\pi,1)$-spaces and so
give natural geometrical model of the classifying spaces for the
braid groups.  Later, V.~I.~Arnold \cite{Arn2}  in this direction 
proved the first results on cohomologies of braids. 
The motivation for his study was a connection (which he discovered)
with the problem of
representing algebraic functions in several variables by
superposing algebraic functions in fewer  variables.
Also, in 1969
V.~I.~Arnold completely described the cohomologies of pure braid
groups \cite{Arn1}. 

In 1969 there appeared the publication of F.~A.~Garside's work \cite{Gar}
where he suggests a new normal form of elements in the braid group
and with its help gives a new solution of the word problem and
also solves the conjugacy problem. In 1968 was published a two-page
note of G.~S.~Makanin \cite{Mak} where he sketches his algorithm
for the solution of the conjugacy problem. The complete publication
of Makanin's work didn't appear (as far as the author is aware).

In 70-ies the study of cohomologies of braids was continued 
independently and by different methods by D.~B.~Fuks \cite{Fuks1} 
who determined 
$\operatorname{mod}~2$ cohomologies, and 
F.~R.~Cohen \cite{CF1},  \cite{CF2},  \cite{CF3} who described 
these homologies with coefficients in $\Bbb Z$ and in
$\Bbb Z/p$ as modules over the Steenrod algebra.

In 1984-85 independently N.~V.~Ivanov \cite{Iv} and J.~McCarthy 
\cite{Mc} proved the
"Tits-alternative" for the mapping class groups of surfaces and 
as a consequence it is true for the braid groups. Namely, 
they proved 
that every subgroup of the mapping class group either contains an 
abelian subgroup of finite index, or
contains a non-abelian free group. 

The question of whether braid groups are linear attracted significant
attention. It was realized that the Burau representation is faithful
for $Br_3$ \cite{Gas}, \cite{MaP}. Then after a long break in 1991
J.~A.~Moody \cite{Mo} proved that  Burau representation is unfaithful
for $n\geq 9$. This bound was improved to $n\geq 6$ by D.~D.~Long 
and M.~Paton \cite{LP} and to $n=5$ by S.~Bigelow \cite{Big1}.
In 1999-2000 there appeared preprints of papers of D.~Krammer 
\cite{Kr1}, \cite{Kr2} and S.~Bigelow \cite{Big2}
who proved that $Br_n$ is linear for all $n$ (using the other
representation).

At the begining of nineties P.~Dehornoy \cite{Deh1}, \cite{Deh2},
\cite{Deh3}  proved
that there exits a left order in braid groups.

Interesting generalizations of braids were introduced in the work of 
E.~Brieskorn \cite{Bri1}. 
The configuration space can also be considered as the orbit space
of the complement of the complexification of the arrangement of
hyperplanes corresponding to the Coxeter group $A_{n-1}=\Sigma_n$.
Generalizing this approach for any finite Coxeter group,
E.~Brieskorn  defined  the so-called generalized braid groups
which are also called Artin groups.

Another way of generalization is to consider braid groups in
3-manifolds, possibly with a boundary. The simplest examples are
braid groups in handlebodies.  A.~B.~Sossinsky \cite{So} was the first
who studied them. Such a group can be interpreted as the
fundamental group of the configuration space of a plane without
$g$ points where $g$ is the genus of the handlebody. 
 The generalized braid
group of type $C$ is isomorphic to the braid group in the solid
torus.

Under the influence of the theory of Vassiliev-Goussarov 
(finite - type) invariants singular braids were 
introduced. The corresponding
algebraic structures are the Baez--Birman monoid \cite{Bae}, 
\cite{Bir2} and
the braid-permutation group by R.~Fenn, R.~Rim\'anyi and 
C.~Rourke
\cite{FRR1}, \cite{FRR2}.
Various properties of these objects were studied \cite{FRZ},
\cite{Zh}, \cite{Ge1}, \cite{Ge2},  \cite{DG}, \cite{Ja}, 
\cite{Cor}, \cite{GM3}.

\section{Definitions and General Properties}
\subsection{Systems of $n$ curves in three-dimensional 
space and braid groups\label{subsec:ini}}

\begin{figure}
\epsfbox{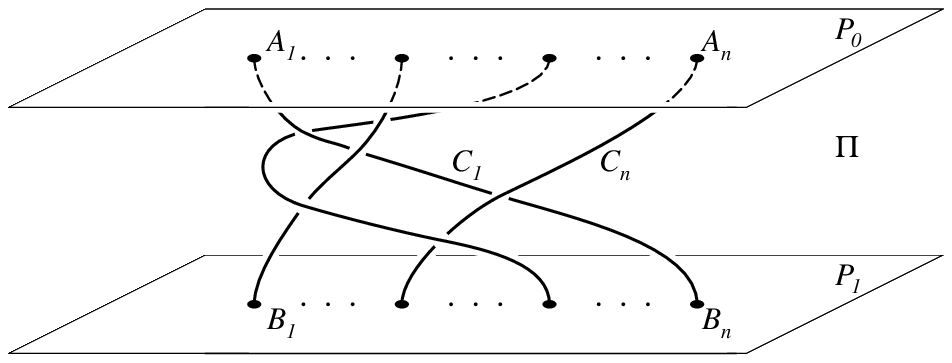}
\caption{} \label{fi:braid}
\end{figure}

First of all as it was already mentioned
braids naturally arise as objects in 3-space. Let us consider two
parallel planes $P_0$ and $P_1$ in $\Bbb R^3$, which contain two
ordered sets of points $A_1$, ..., $A_n\in P_0$ and $B_1$, ...,
$B_n\in P_1$. These points are lying on parallel lines $L_A$ and
$L_B$ respectively. The space between the planes $P_0$ and $P_1$
we denote by $\Pi$. Suppose that the point $B_i$ is lying under
the point $A_i$, as a result of the orthogonal projection of the plane
$P_0$ onto the plane $P_1$. Let us connect the set of points $A_1$, ...,
$A_n$ with the set of points $B_1$, ..., $B_n$ by simple non-intersecting
curves
$C_1$, ..., $C_n$ lying in the space $\Pi$ and such that each
curve meets only once each parallel plane $P_t$
 lying in the space
$\Pi$ (see Figure~\ref{fi:braid}). 
This object is called a {\it
braid} and the curves are called the {\it strings} of a braid.
Usually braids are depicted by projections on the plane passing
through the lines $L_A$ and $L_B$. This projection is supposed to
be in general position so that there is only finite number of
double points of intersection which are lying on pairwise
different levels and intersections are transversal. The simplest
braid $\sigma_i$ (Fig.~\ref{fi:sigma}) corresponds to the transposition
$(i,i+1)$.

\begin{figure}
\epsfbox{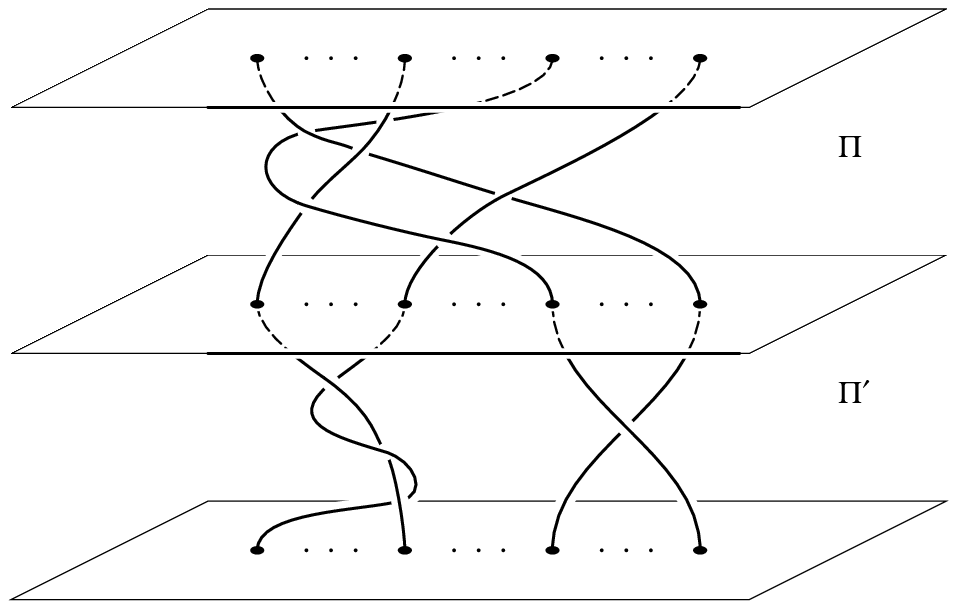}
\caption{}\label{fi:can}
\end{figure}

 Let us introduce the following equivalence relation on the
set of all braids with $n$ strings and with fixed $P_0$, $P_1$,
$A_i$ and $B_i$. It is defined by homeomorphisms $h:\Pi\to \Pi$,
identical on $P_0\cup P_1$ and such that $h(P_t)=P_t$. Braids
$\beta$ and $\beta^\prime$ are equivalent if there exists such a
homeomorphism $h$ that $h(\beta)=\beta^\prime$.
\begin{figure}
\vskip 0.2cm
\begin{picture}(0,130)(0,-10) \thicklines
\put(-100,100){\line(0,-1){100}} \put(-50,100){\line(0,-1){100}}
\put(-25,100){\line(1,-2){50}} \put(25,100){\line(-1,-2){20}}
\put(-25,0){\line(1,2){20}} \put(50,100){\line(0,-1){100}}
\put(100,100){\line(0,-1){100}}
\put(-100,110){\makebox(0,0)[cc]{$1$}}
\put(-50,110){\makebox(0,0)[cc]{$i-1$}}
\put(-25,110){\makebox(0,0)[cc]{$i$}}
\put(25,110){\makebox(0,0)[cc]{$i+1$}}
\put(50,110){\makebox(0,0)[cc]{$i+2$}}
\put(100,110){\makebox(0,0)[cc]{$n$}}
\put(-75,50){\makebox(0,0)[cc]{.\quad.\quad.}}
\put(75,50){\makebox(0,0)[cc]{.\quad.\quad.}}
\end{picture}
\caption{}\label{fi:sigma}
\end{figure} 
On the set $Br_n$ of equivalence classes under the considered
relation the structure of a group introduces as follows. We put a
copy $\Pi^\prime$ of the domain $\Pi$ under the $\Pi$ in such a
way that $P_0^\prime$ coincides with $P_1$ and each $A_i$
coincides with $B_i$ and we glue braids $\beta$ and
$\beta^\prime$. This gluing gives a composition of braids $\beta
\beta^\prime$ (Fig.~\ref{fi:can}). 
Unit element is the equivalence class
containing a braid of $n$ parallel intervals, the braid
$\beta^{-1}$ inverse to $\beta$ is defined by reflection of
$\beta$ with respect to the plane $P_{1/2}$. 
A string $C_i$ of a braid $\beta$ connects the point $A_i$ with
the pont $B_{k_i}$ defining the permutation $S^\beta $. If this
permutation is identical then the braid $\beta$ is called {\it
pure}.
The map $\beta\to
S^\beta$ defines an epimorphism $\tau_n$ of the braid group $Br_n$
on the permutation group $\Sigma_n$with the kernel consisting of
all pure braids: 
\begin{equation}
1\rightarrow P_n\rightarrow Br_n \buildrel
\tau_n \over \longrightarrow \Sigma_n\rightarrow 1 . \label{eq:exact}
\end{equation}

The 
following presentation of the braid group $Br_n$ with generators $\sigma_i$, 
$i=1, ..., n-1$ and two types of relations: 
\begin{equation}
 \begin{cases} \sigma_i \sigma_j &=\sigma_j \, \sigma_i, \ \
\text{if} \ \ |i-j|
>1,
\\ \sigma_i \sigma_{i+1} \sigma_i &= \sigma_{i+1} \sigma_i \sigma_{i+1}
\end{cases} \label{eq:brelations}
\end{equation}
is the algebraic expression of the fact that any isotopy of braids
can be  broken down
into ``elementary moves'' of two types that correspond to two types of
relations.

If we add a vertical interval to the system of curves on 
Figure~\ref{fi:braid} we can get a canonical inclusion $j_n$ of the group
$Br_n$ into the group $Br_{n+1}$
$$j_n:Br_{n} \to Br_{n+1}.$$
If the symmetric group $\Sigma_n$ is given by its canonical presentation
with generators $s_i$, 
$i=1, ..., n-1$ and relations: 
\begin{equation}
 \begin{cases} s_i s_j &=s_j \, s_i, \ \
\text{if} \ \ |i-j|
>1,
\\ s_i s_{i+1} s_i &= s_{i+1} s_i s_{i+1}\\
s_i^2=1,
\end{cases} \label{eq:srelations}
\end{equation}
then the homomorphism $\tau_n$ is given by the formula
$$\tau_n(\sigma_i) = s_i, \ \ \ i= 1, \dots , n-1.$$

It is possible to consider braids as classes of equivalence of
\emph{braid diagrams} which are generic projections of three dimensional
braids on a plane. The classes of equivalence are defined
by the \emph{Reidemeister moves} depicted at Figure~\ref{fi:R23}.
\begin{figure}
\begin{picture}(0,110)(0,-15)
\put(-120,80){\line(1,-2){20}} \put(-100,80){\line(-1,-2){8}}
\put(-70,80){\line(0,-1){80}} \put(-50,80){\line(0,-1){80}}
\put(0,80){\line(1,-2){40}} \put(20,80){\line(-1,-2){8}}
\put(40,80){\line(-1,-2){18}} \put(80,80){\line(1,-2){40}}
\put(100,80){\line(1,-2){20}} \put(120,80){\line(-1,-2){8}}
\put(102,44){\line(1,2){6}} \put(-120,40){\line(1,2){8}}
\put(-120,40){\line(1,-2){8}} \put(-100,40){\line(-1,-2){20}}
\put(0,40){\line(1,2){8}} \put(0,40){\line(1,-2){20}}
\put(120,40){\line(-1,-2){8}} \put(12,24){\line(1,2){6}}
\put(-100,0){\line(-1,2){8}} \put(0,0){\line(1,2){8}}
\put(80,0){\line(1,2){18}} \put(100,0){\line(1,2){8}}
\put(-85,40){\makebox(0,0)[cc]{$\leftrightarrow$}}
\put(60,40){\makebox(0,0)[cc]{$\leftrightarrow$}}
\end{picture}
\caption{}\label{fi:R23}
\end{figure} 

\subsection{Braid groups and configuration spaces\label{subsection:braco}}

\begin{figure}
\epsfbox{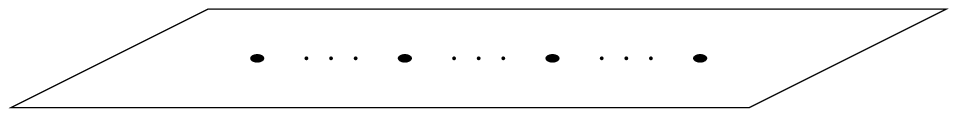}
\caption{} \label{fi:conf}
\end{figure}

If we  look at the Figure \ref{fi:braid}, then this picture
can be interpreted as a graph of 
 a loop in the {\it configuration space} of $n$ points 
on a plane, that is the space of unordered sets of $n$ points 
on a plane, see Figure~\ref{fi:conf}.
So, it is possible to interpret the braid group as 
the fundamental group of the configuration space.
Formally it is done as follows.
The symmetric group $\Sigma_m$ acts on the
Cartesian power $(\R^2)^m$ of the space $\R^2$: 
\begin{equation}
w(y_1,...,y_m) =
(y_{w^{-1}(1)},...,y_{w^{-1}(m)}), \ \ w \in \Sigma_m.
\label{eq:permut}
\end{equation}
Denote
by $F(\R^2,m)$ the space of $m$-tuples of pairwise different points
in $\R^2$: 
$$F(\R^2,m)=\{(p_1,...,p_m)\in (\R^2)^m:p_i\not=p_j \ \text {for} \
i\not=j \}.$$ 
This is the space of regular points of our action.
 We call the orbit space of this action  $B(\R^2,m) = F(\R^2,m)/\Sigma_m$
the {\it configuration
space of} $n$ {\it points on a plane}. The braid group
$Br_m$ is the fundamental group of configuration space
$$Br_m= \pi_1(B(\R^2,m)). $$
The pure braid group $P_m$ is the 
is the fundamental group of  the space  $F(\R^2,m)$. 
The covering 
 $$p:F(\R^2,m)\rightarrow B(\R^2,m)$$ 
defines the exact sequence:
\begin{equation}
1 \rightarrow\pi_1(F(\R^2,m))\buildrel p_*
\over\rightarrow \pi_1(B(\R^2,m)) \rightarrow \Sigma_n
\rightarrow 1 ,\label{eq:pi1conf}
\end{equation}
which is equivalent to 
sequence (\ref{eq:exact}).

  It can be used for proving the canonical presentation of the
braid group (\ref{eq:brelations}) as it is done, for example in 
the book of J.~Birman \cite{Bir1}.

Such considerations were made by R.~Fox and L.~Neuwirth
\cite{FoN}. 

\subsection{Braid groups as automorphism groups of free groups and 
the word problem\label{subsec:aut}}

Another important approach to the braid group bases on the fact
that this group may be considered as a subgroup of the
automorphism group of a free group.

Let $F_n$ be the free group of rank $n$ with the set of generators
$\{x_1, \dots, x_n \}$. Assume further that $\operatorname{Aut}
F_n$ is the automorphism group of $F_n$. 

We have the standard
inclusions of the symmetric group $\Sigma_n$ and the braid group
$Br_n$ into $\operatorname{Aut} F_n$. For the braid group it may be
described as follows. 
Let $\overline{\sigma}_i\in \operatorname{Aut} F_n, i=1,2,
\dots, n-1,$ be given by the formula which describes its action on
generators:
\begin{equation} \begin{cases} 
x_i &\mapsto x_{i+1},
\\ x_{i+1} &\mapsto x_{i+1}^{-1}x_ix_{i+1}, \\
x_j &\mapsto x_j, j\not=i,i+1. 
\end{cases} \label{eq:autf}
\end{equation}
Let us define a map  $\nu$ of the generators $\sigma_i$, 
$i=1, \dots, n-1$ of the braid group $Br_n$ to these automorphisms:
\begin{equation} 
\nu(\sigma_i)= \overline{\sigma}_i.
\label{eq:incaut}
\end{equation}
\begin{Theorem} Formulas  (\ref{eq:incaut}) define correctly a homomorphism
$$\nu: Br_n
 \rightarrow \operatorname{Aut} F_n .$$
which is a monomorphism.
\label{Theorem:theaut}
\end{Theorem} 
Theorem \ref{Theorem:theaut} gives a solution of the word problem for 
the braid groups. It was done first by E.~Artin \cite{Art1}. 

The free group
$F_n$ is a fundamental group of a disc $D_n$ without $n$ points
and the generator $x_i$ corresponds to a loop going around the
$i$-th point. The braid group $Br_n$ is the mapping class group of
a disc $D_n$ with its boundary fixed \cite{Bir1} and so it acts on the
fundamental group of $D_n$.
This action is described by the formulas (\ref{eq:autf}) where
$x_i$ correspond to the canonical loops on $D_n$ which form the
generators of the fundamental group. Geometrically
this action is depicted in the Figure~\ref{fi:mapcl}.

\begin{figure}
\input mapcl
\caption{}\label{fi:mapcl}
\end{figure} 

\subsection{Commutator subgroup and other presentations\label{subsec:commu}}

Let us define a homomorphism from braid group to
integers by taking the sum of exponents of the entries of the
generators $\sigma_i$ in the expression of any element of the
group through these canonical generators: 
$$\operatorname{deg}:
Br_n \rightarrow \Bbb Z, \ \operatorname{deg}(b) = \sum_j m_j, \
\text{where} \ b=(\sigma_{i_1} )^{m_1} \dots
(\sigma_{i_k})^{m_k}.$$ 
\begin{Proposition} The homomorphism 
$$\operatorname{deg}:
Br_n \rightarrow \Bbb Z$$
gives the abelianization of the braid group and the commutator subgroup 
$Br_n^\prime$
is characterized by the condition
$$b\in Br_n^\prime \ \ {\it if} \ {\rm and} \  {\it only} \ {\it if} \ \ 
\operatorname{deg}(b) =0.$$
\end{Proposition}
\begin{proof} Let $a:Br_n\to A$ be a homomorphism to any other
abelian group $A$, then from the relations (\ref{eq:brelations})
we have:
$$a(\sigma_i) a(\sigma_{i+1}) a(\sigma_i) = a(\sigma_{i+1})
 a(\sigma_i) a(\sigma_{i+1}).$$
the commutativity of $A$ gives that 
$a(\sigma_{i+1}) =a(\sigma_{i}) .$ This means that the homomorphism
$\operatorname{deg}$ is universal. 
\end{proof}

Of course, there exist another presentations of the braid group.
Let $$\sigma = \sigma_1 \sigma_{2} \dots \sigma_{n-1},$$ then the
group $Br_n$ is generated by $\sigma_1$ and $\sigma$ because
$$\sigma_{i+1} =\sigma^i \sigma_1 \sigma^{-i}, \quad i =1, \dots
{n-2}.$$
The relations for the generators $\sigma_1$ and $\sigma$ are the 
following
\begin{equation}
 \begin{cases}
\sigma_1 \sigma^i \sigma_1 \sigma^{-i} &= 
\sigma^i \sigma_1 \sigma^{-i} \sigma_1 \ \  \text{for} \ \
2 \leq i\leq {n / 2}, \\
\sigma^n &= (\sigma \sigma_1)^{n-1}.
\end{cases} \label{eq:2relations}
\end{equation}
This was observed by Artin in the initial paper \cite{Art1}.

An interesting series of presentations was given by V.~Sergiescu
\cite{Ser}. For every planar graph he constructed a presentation of the
group $Br_n$, where $n$ is the number of vertices of the graph,
with generators corresponding to edges and relations reflecting
the geometry of the graph. Artin's presentation in this context
corresponds to the graph consisting of the interval from 1 to $n$
with the natural numbers (from 1 to $n$) as vertices and with
segments between them as edges. For gereralizations of braids graph
presentations of these type were considered by P.~Bellingeri and 
V.~Vershinin \cite{Bel, BelV}

\begin{figure}
\epsfbox{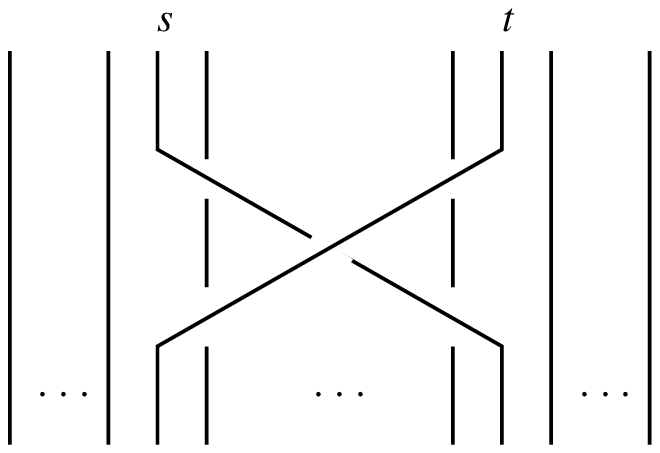}
\caption{}
\label{fi:sbige2}
\end{figure}

J.~S.~Birman;  K.~H.~Ko; S.~J.~Lee \cite{BKL} introduced the presentation 
with the generators $a_{ts}$ with $1 \leq s<t\leq n$
and relations

\begin{equation} \begin{cases}
a_{ts}a_{rq}&=a_{rq}a_{ts} \ \ {\rm for} \ \ (t-r)(t-q)(s-r)(s-q)>0,\\ 
a_{ts}a_{sr} &=a_{tr}a_{ts}=a_{sr}a_{tr}  \ \ {\rm for} \ \ 
1\leq r<s<t\leq n .
\end{cases}\label{eq:rebkl}
\end{equation}
The generators $a_{ts}$ are expressed by the canonical generators 
$\sigma_i$
in the following form:
 \begin{equation} a_{ts}=(\sigma_{t-1}\sigma_{t-2}\cdots\sigma_{s+1})\sigma_s
(\sigma^{-1}_{s+1}\cdots\sigma^{-1}_{t-2}\sigma^{-1}_{t-1})  \ \ 
{\rm for} \ \ 1\leq s<t\leq n.
\label{eq:ats}
\end{equation} 
Geometrically the generators $a_{s,t}$ are depicted in Figure~\ref{fi:sbige2}.

The set of generators for braid groups was even enlarged 
in the work of Jean Michel \cite{Mi} as follows.
Let $\vert \, \, \vert : \Sigma_n\to \Z$
be the length function on the symmetric group with respect to the 
generators $s_i$: for $x\in \Sigma_n$,
$\vert  x \vert $ is the smallest natural number $k$ such that $x$ 
is a product of $k$ elements of the set 
$\{s_1,...,s_{n-1}\}$.  
It is known (\cite{Bo}, Sect.~1, Ex.~13(b)) that two minimal expressions
for an element of $\Sigma_n$ are equivalent by using only of the 
relations (\ref{eq:brelations}). This implies that
the canonical projection $\tau_n :Br_n\to \Sigma_n$ has a  unique  
set-theoretic section
$r:\Sigma_n\to Br_n$ 
such that $r(s_i)=\sigma_i$ for $i=1,...,n-1$
and $r(xy)= r(x)\, r(y)$ whenever 
$\vert  xy\vert =\vert  x\vert +\vert  y\vert$. Then the group 
$Br_n$ admits a
presentation by generators
$\{ r(x) \, \vert\,  x\in \Sigma_n\}$ and relations $r(xy)= r(x) \,r(y)$  
for all $x,y\in \Sigma_n$ such that
 $\vert  xy\vert =\vert  x\vert +\vert  y\vert$.

\subsection{Presentation of the pure braid group and Markov normal 
form\label{subsec:pure}}

Let $f(y_1,..., y_m$) be a word with (possibly empty) entries of
$y_i^\epsilon$, where $y_i$ are some letters and $\epsilon$ may be
$\pm 1$. If $y_i$ are elements of a group $G$ then $f(y_1,...,
y_m$) will be considered as the corresponding element of $G$.

Let us define the elements $s_{i,j}$, $1\leq i<j\leq m$, of the
braid group $Br_m$ by the formula:
$$s_{i,j}=\sigma_{j-1}...\sigma_{i+1}\sigma_{i}^2\sigma_{i+1}^{-1}...
\sigma_{j-1}^{-1}.$$ 
These elements satisfy the following Burau
relations (\cite{Bu1}, \cite{Mar2}, see also \cite{Li1}): 
\begin{multline}
\begin{cases}
s_{i,j}s_{k,l}&=s_{k,l}s_{i,j}
\ \text {for} \ i<j<k<l \ \text {and} \ i<k<l<j, \\
s_{i,j}s_{i,k}s_{j,k}&=s_{i,k}s_{j,k}s_{i,j} \ \text {for} \
i<j<k, \\
s_{i,k}s_{j,k}s_{i,j}&=s_{j,k}s_{i,j}s_{i,k} \ \text
{for} \ i<j<k, \\
s_{i,k}s_{j,k}s_{j,l}s_{j,k}^{-1}&=s_{j,k}s_{j,l}s_{j,k}^{-1}s_{i,k}
\ \text {for} \ i<j<k<l.\\
\end{cases}
\label{eq:burau}
\end{multline}
W.~Burau and later A.~A.~Markov proved that the elements
$s_{i,j}$ with the relations (\ref{eq:burau}) give a~presentation of the pure
braid group $P_m$ \cite{Mar2}. The following formula is a consequence
of the Burau relations and also belongs to A.~A.~Markov:
\begin{equation}
[s_{i,l},s_{j,k}^\epsilon]=f(s_{1,l},...,s_{l-1,l}), \ \
\epsilon=\pm 1, \ \ k<l.\label{eq:commut}
\end{equation}
 Let us define the elements
$\sigma_{k,l}$, $1\leq k\leq l\leq m$ by the formulas
$$\sigma_{k,k}=e,$$ 
$$\sigma_{k,l}=\sigma_{k}^{-1}...\sigma_{l-1}^{-1}.$$
Let $P^k_m$ be the subgroup of $P_m$ generated by the elements
$s_{i,j}$ with $k<j$. 
\begin{Theorem} {\rm(A.~A.~Markov)}
(i) Every element of the group $Br_{m}$ can be uniquely written in
the form
\begin{equation}
f_{m}(s_{1,m},...,s_{m-1,m})...f_j(s_{1,j},...,s_{j-1,j})...f_2(s_{1,2})
\sigma_{i_{m},m}...\sigma_{i_{j},j}...\sigma_{i_{2},2}.\label{eq:markform}
\end{equation}

(ii) The factor group $P^k_m/P^{k-1}_m$ is the free group on free
generators $s_{i,{k+1}}$, $1\leq i\leq k$.
\label{Theorem:mnf}
\end{Theorem}
The form (\ref{eq:markform}) is called the \emph{Markov normal form},
it also gives the solution 
of the word problem for the braid groups.

\section{Garside Normal Form, Center and Conjugacy Problem\label{sec:gar}}

Essential role in Garside work \cite{Gar} plays the monoid of 
\emph{positive} braids $Br_n^+$, that is the monoide which has a presentation
with generators  $\sigma_i$, 
$i=1, ..., n$ and relations (\ref{eq:brelations}). In other words each
element of this monoide can be represented as a word on the elements 
$\sigma_i$, $i=1, ..., n$ with no entrances of  $\sigma_i^{-1}$.
Two positive words $A$ and $B$ in the
alphabet $\{\sigma_i$, $(i=1,\dots,n-1) \}$
will be said to be {\it positively equal} if they are equal as elements
of $Br_n^+$. In this case we shall write $A\doteq B$.

First of all Garside proves the following statement.
\begin{Proposition} In $Br_n^+$ for $i,k = 1, ..., n-1 $, given
$\sigma_i A \doteq \sigma_k B$,
it follows that
\begin{equation*} \text{if} \ \  k = i, \ \text{then} \ \ A \doteq B,
\end{equation*}
\begin{equation*} \text{if} \ \  |k - i| = 1,
\ \text{then} \ \ A \doteq \sigma_k\sigma_i Z, \ B\doteq \sigma_i\sigma_k Z
\ \ \text{for some} \ \ Z,
\end{equation*}
\begin{equation*} \text{if} \ \ |k - i| \geq 2, \ \text{then} \ \ A \doteq
\sigma_k Z, \  B\doteq \sigma_i Z \ \ \text{for some} \ \ Z.
\end{equation*}
The same is true for the right multiples of $\sigma_i$.
\label{Proposition:div}
\end{Proposition} 
\begin{Corollary} If $A\doteq P$, $B\doteq Q$, $AXB \doteq PYQ$, ($L(A)\geq 0$,
$L(B) \geq 0$), then $X\doteq Y$. That is, monoid $Br_n^+$ is left and right
cancellative.
\end{Corollary}
Garside's  {\it fundamental word} $\Delta$ in the braid 
group $Br_{n+1}$ is defined by the formula:
$$\Delta = \sigma_1 \dots \sigma_n \sigma_1 \dots \sigma_{n-1} \dots  
\sigma_1 \sigma_2 \sigma_1.$$
If we use Garside's notation $\Pi_t\equiv \sigma_1\dots \sigma_t$, then
$\Delta \equiv \Pi_{n-1} \dots \Pi_1$.

For a positive word $W$ in $\sigma_i$, $i=1, ..., n$ we say that $\Delta$
is a \emph{factor} of $W$ or simply $W$ contains $\Delta$, if 
$W \doteq A \Delta B$ with $A$ and $B$ being
arbitrary positive words, probably empty.
If $W$ does not contain $\Delta$ we shall say $W$ is \emph{prime to} 
$\Delta$.

Garside's transformation of words $\mathcal R$ is defined by the formula
\begin{equation*}
\mathcal R(\sigma_i) \equiv \sigma_{n-i}.
\end{equation*} 
This gives  the
automorphism of $Br_n$ and the positive braid monoid $Br_n^+$.
\begin{Proposition} In $Br_n$ 
$$\sigma_i\Delta \doteq\Delta \mathcal R(\sigma_i).$$
\end{Proposition}
Geometrically this commutation is  shown on Figure~\ref{fi:deltasi} 
($\Delta\sigma_3 = \sigma_1\Delta $).
\begin{figure}
\begin{picture}(0,230)(0,-20)
\thicklines
\put(-120,200){\line(1,-1){12}} 
\put(-90,200){\line(-1,-1){30}}
\put(-60,200){\line(0,-1){30}} 
\put(-30,200){\line(0,-1){60}}
\put(50,200){\line(1,-1){12}} 
\put(80,200){\line(-1,-1){30}}
\put(110,200){\line(0,-1){60}} 
\put(140,200){\line(0,-1){90}}
\put(-104,182){\line(1,-1){24}} 
\put(-60,170){\line(-1,-1){30}}
\put(-120,170){\line(0,-1){60}}
\put(50,170){\line(1,-1){12}}
\put(80,170){\line(-1,1){12}} 
\put(80,170){\line(-1,-1){30}} \put(140,170){\line(0,-1){60}}
\put(66,152){\line(1,-1){24}}
\put(-74,152){\line(1,-1){24}}
\put(-90,140){\line(0,-1){30}}
\put(-30,140){\line(-1,-1){30}}
\put(50,140){\line(0,-1){60}} 
\put(110,140){\line(-1,-1){30}}
\put(96,122){\line(1,-1){24}}
\put(-90,110){\line(-1,-1){30}}
 \put(-120,110){\line(1,-1){12}}
\put(80,110){\line(0,-1){30}} 
\put(140,110){\line(-1,-1){30}}
\put(-30,110){\line(-1,1){12}} 
\put(-30,110){\line(0,-1){90}} 
\put(-102,92){\line(1,-1){24}} 
\put(-120,80){\line(0,-1){30}}
\put(-60,80){\line(-1,-1){60}} 
\put(-60,110){\line(0,-1){30}}
\put(50,80){\line(1,-1){12}} 
\put(80,80){\line(-1,-1){30}} 
\put(140,80){\line(-1,1){12}}
\put(140,80){\line(0,-1){90}}
\put(-72,62){\line(1,-1){12}} 
 \put(68,62){\line(1,-1){24}} 

\put(50,50){\line(0,-1){30}} 
\put(110,80){\line(0,-1){30}}
\put(-60,20){\line(1,-1){12}}
\put(-42,2){\line(1,-1){12}}
\put(50,20){\line(1,-1){12}} 
\put(-120,50){\line(1,-1){12}} 
\put(-90,20){\line(0,-1){30}} 
\put(-120,20){\line(0,-1){30}} 
\put(-30,20){\line(-1,-1){30}}
\put(110,20){\line(0,-1){30}} 
\put(110,50){\line(-1,-1){60}}
\put(-60,50){\line(0,-1){30}} 
\put(98,32){\line(1,-1){12}}
\put(68,2){\line(1,-1){12}}
\put(-102,32){\line(1,-1){12}}
\put(10,90){\makebox(0,0)[cc]{$-$}}
\put(10,85){\makebox(0,0)[cc]{$-$}}
\end{picture}
\caption{}\label{fi:deltasi}
\end{figure}

\begin{Proposition} If $W$ is an arbitrary positive word in $Br_n^+$ such that either
\begin{equation*}
W\doteq \sigma_1 A_1\doteq\sigma_2 A_2\doteq\dots\doteq\sigma_{n-1}A_{n-1},
\end{equation*}
or
\begin{equation*}
W\doteq B_1\sigma_1 \doteq B_2\sigma_2 \doteq\dots\doteq
B_{n-1}\sigma_{n-1},
\end{equation*}
then $W\doteq\Delta Z$ for some $Z$.
\end{Proposition}
\begin{Proposition} The canonical homomorphism
$$Br_n^+ \to Br_n$$
is a monomorphism.
\label{Propositionro:inj}
\end{Proposition}
Among positive words on the alphabet  $\{\sigma_1 \dots \sigma_n\}$ let us 
introduce a lexicographical ordering with the condition that  
$\sigma_1 < \sigma_2 < \dots < \sigma_n $. For a positive word $W$ the 
\emph{base} of $W$ is the smallest positive word which is positively equal
to $W$. The base is uniquely determined. If a positive word $A$ is prime 
to $\Delta$, then for the base of $A$ the notation $\overline{A}$  will
be used.
\begin{Theorem} {\rm(F.~A.~Garside)}
 Every word $W$ in $Br_{n+1}$ can be uniquely written in
the form $\Delta^m \overline{A}$, where $m$ is an 
 integer.
\label{Theorem:garnf}
\end{Theorem}
The form of a word $W$ established in this theorem we call the
\emph{Garside left normal form} and the 
index $m$ we call the \emph{power} of $W$. The same way the
\emph{Garside right normal form} is defined and the corresponding
variant of Theorem~\ref{Theorem:garnf} is true. The
Garside normal form also gives a solution to the word problem in the braid
group.
\begin{Theorem} {\rm (F.~A.~Garside)}
The necessary and sufficient condition that two words in  
$Br_{n+1}$ are equal is that their Garside normal forms are identical.
\label{Theorem:gws}
\end{Theorem}
Garside normal form for the braid groups was precised in the subsequent 
works of S.~I.~Adyan \cite{Ad},
 W.~Thurston \cite{E_Th}, E.~El-Rifai and H.~R.~Morton \cite{EM}. 
Namely, there was introduced 
the \emph{left-greedy form} (in the terminology of
W.~Thurston \cite{E_Th}) 
\begin{equation*}
\Delta^t A_1 \dots A_k,
\end{equation*}
where $A_i$ are the successive
possible longest \emph{fragments of the word} $\Delta$ (in the terminology 
of S.~I.~Adyan \cite{Ad}) or \emph{positive permutation braids} (in the 
terminology of E.~El-Rifai and H.~R.~Morton \cite{EM}). Certainly, the same 
way the \emph{right-greedy form} is defined. With the help of this form
it proved that the braid group is biautomatic.

The center of the braid group was firstly found by W.-L.~Chow \cite{Ch}. 
Namely, as it follows from the presentation of braid groups with two
generators $\sigma_1$ and $\sigma$ and relations (\ref{eq:2relations})
given in the subsection  \ref{subsec:ini}
the element $\sigma^n$ commutes with $\sigma \sigma_1$ and so with 
$\sigma_1$. Chow proved that it generates the center.
 Garside normal form gives an elegant proof of the following theorem.
\begin{Theorem} 
(i) When $n=1$ the center of the group $Br_{n+1}$ is generated by $\Delta$.
\par
(ii) When $n>1$ the center of the group $Br_{n+1}$ is generated by 
$\Delta^2$.
\label{Theorem:centr}
\end{Theorem}
Let $\alpha$ be a positive word such that $\Delta \doteq \alpha X$, where
$X$ is an arbitrary positive word, probably empty. For any word $W$ in
$Br_{n+1}$, the word $\alpha^{-1}W\alpha$, reduced to Garside normal form
is called an $\alpha$-\emph{transformation of} $W$. 
 
For any word $W$ in $Br_{n+1}$ with the Garside
normal form $\Delta^m \overline{A} \equiv W_1$ consider the following
chains of $\alpha$-transformations: take all the $\alpha$-transformations 
of $W_1$ and let those which are of power $\geq m$ and which are distinct
from each other be $W_2$, $W_3$, $\dots$, $W_t$. Now repeat the process 
for each of the words $W_2$, $W_3$, $\dots$, $W_t$ in turn, denoting
successively by $W_{t+1}$, $W_{t+2}$, $\dots$, any new words occurring, the
condition being always that each new word must be of power $\geq m$.
Continue to repeat the process for every new distinct word arising, as
the sequence  $W_1$, $W_2$,  $W_{t+2}$, $\dots$, expands.
\begin{Proposition} 
The set $W_1$, $W_2$,  $W_{t+2}$, $\dots$, is finite.
\end{Proposition}
 Suppose that in the set  $W_1$, $W_2$,  $W_{t+2}$, $\dots$, the highest 
power reached is $s$ and that the words of power $s$ form the subset
 $V_1$, $V_2$,  $\dots$. Then this set $V_1$, $V_2$,  $\dots$ is
called the \emph{summit set} of $W$.
\begin{Theorem} {\rm(F.~A.~Garside)}
 Two elements $A$ and $B$ of the group $Br_{n+1}$ are conjugate if and 
only if their summit sets are identical.
\label{Theorem:cong}
\end{Theorem}

J.~S.~Birman, K.~H.~Ko and S.~J.~Lee considered
the word $\delta=a_{n(n-1)}\cdots
a_{32}a_{21}=\sigma_{n-1}\cdots\sigma_2\sigma_1$, 
 as  a fundamental in their system of generators  and proved that 
every element in $Br_n$ has a representative $W=\delta^jA_1A_2\cdots
A_k$ with positive $A_i$ in a unique way in some sense. Based on this form
they gave an algorithm for the word problem in $B_n$
which runs in time $(nح^2)$ for a given word of length $m$.

\section{Ordering of Braids\label{sec:dehor}}

A group $G$ is said \emph{totally} (or \emph{linearly}) \emph{left}
(correspondingly \emph{right}) \emph{ordered} 
if it has a total order $<$ invariant by left (right) multiplication,
i.e. if $a<b$, then $ca <  cb$ for any $ c \in G$. If this order 
is also invariant by right (left) multiplication, then the group
$G$ is called \emph{ordered}.

For any left ordered group $G$ denote by $P$ the set of positive 
elements $\{x\in G: x > 1\}$, then the set of negative elements is 
defined by the formula: $P^{-1} = \{x\in G: x\in P\}$. The total character 
of an order on $G$
is expressed by the partition
$$G= P \coprod \{1\} \coprod P^{-1}.$$
The invariance of multiplication is expressed by the inclusion 
$P^2 \subset P$, where $P^2$ is formed by products of couples of 
elements of $P$. Conversely, if there exist a subset $P$ of  a 
group $G$ with the properties:
$$G= P \coprod \{1\} \coprod P^{-1}, \quad P^2 \subset P, $$
then $G$ is left ordered by the order defined by: 
$x < y$ if and only if $x^{-1}y \in P$. A group $G$ then is ordered 
if and only if $xPx^{-1} \subset P$ for all $x \in G$.

Let $i\in \{1, \dots, n\}$ and a word $w$ on the alphabet 
$\{ \sigma_1, \dots, \sigma_n\}$ is expressed in the form
$$w_0 \sigma_i w_1 \sigma_i \dots \sigma_i w_r,$$ 
where the
subwords $w_0, \dots, w_r$ are the words on the letters 
$\sigma_j^{\pm 1}$ with $j> i$. Then such a word is called 
$\sigma_i$-\emph{positive}. This means that all entries
of $\sigma_i^{\pm 1}$ in the word $w$ with $i$ minimal
must be positive. If all such entries are negative then
a word $w$ is called $\sigma_i$-\emph{negative}.
A braid of $Br_{n+1}$ is called 
$\sigma_i$-\emph{positive} ($\sigma_i$-\emph{negative}) if there exists 
its expression as
a word on
the standard generators which is $\sigma_i$-\emph{positive} 
($\sigma_i$-\emph{negative}). A braid is called 
$\sigma$-\emph{positive} ($\sigma$-\emph{negative})  if it exists
a number $i$, such that it is $\sigma_i$-\emph{positive} 
($\sigma_i$-\emph{negative}). 
\begin{Theorem} {\rm(P.~Dehornoy)}
 Every braid in $Br_{n+1}$ different from 1 is either 
$\sigma$-positive or $\sigma$-negative.
\label{Theorem:deh}
\end{Theorem}
\begin{Corollary} For all $n$ the braid group  $Br_{n+1}$ is left
ordered.
\end{Corollary}

\section {Representations\label{sec:rep}}
\subsection{Burau representation\label{subsec:bur}}

Let us map the generators of the braid group $Br_n$ to
the following elements of the group $GL_n\Bbb Z[t,t^{-1}]$
\begin{equation}
\sigma_i \mapsto \left( \begin{matrix} E_{i-1} & 0 & 0 & 0 \\ 0 & 1-t &
t & 0 \\ 0 & 1   & 0 & 0 \\ 0 & 0   & 0 & E_{n-i-1}
\end{matrix}
\right), \label{eq:bur}
\end{equation}
where $E_i$ is the unit $i\times i$ matrix.
The formula (\ref{eq:bur}) defines correctly
the representation of the braid group in $GL_n\Bbb
Z[t,t^{-1}]$: $$r: Br_n\to GL_n\Bbb Z[t,t^{-1}],$$
which is called {\it Burau representation} \cite{Bu3}.
\begin{Theorem}  Burau representation is faithful for 
$n = 3$. \label{Theorem:bufai}
\end{Theorem}
\begin{Theorem} ({\rm J.~A.~Moody;  D.~D.~Long and M.~Paton;  
S.~Bigelow}) 
Burau representation is unfaithful
for $n\geq 5$.  \label{Theorem:bunf}
\end{Theorem}
The case $n=4$ remains open.

\subsection{Lawrence-Krammer representation\label{subsec:lakre}}

Consider the ring $K= \Z[q^{\pm 1}, t^{\pm 1}]$ of Laurent polynomials 
in two variables $q,t$, and the
 free $K$-module $$ V= \bigoplus_{1 \le i<j \le n} K\,x_{i,j}. $$ 
For $k \in \{1, 2, \dots, n-1\}$, define the action of the braid 
generators $\sigma_k$ on the basis of $V$ by the formula:
\begin{equation} \sigma_k(x_{i,j}) = \begin{cases}
 x_{i,j},  &k<i-1 \ {\rm or} \ j<k ;\\ 
 x_{i-1,j} +(1-q) x_{i,j}, &k=i-1;\\ 
tq(q-1) x_{i,i+1} +q x_{i+1,j},&k=i<j-1;\\ 
tq^2 x_{i,j}; &k=i=j-1;\\ 
x_{i,j}+tq^{k-i}(q-1)^2 x_{k,k+1}, &i<k<j-1;\\
x_{i,j-1}+tq^{j-i}(q-1) x_{j-1,j}, &k=j-1;\\
 (1-q)x_{i,j} +q x_{i,j+1}, &k=j .\end{cases}\label{eq:lakra}
\end{equation}
Direct computation shows that this defines a representation
$$\rho_n:Br_n\to GL(V),$$
which was firstly defined by R.~Lawrence \cite{Law} in topological
terms and in the explicit form (\ref{eq:lakra}) by D.~Krammer \cite{Kr2}.
\begin{Theorem} {\rm (S. Bigelow \cite{Big2}, D.~Krammer \cite{Kr2})}  
The representation $$\rho_n : Br_n\to GL \, (V)$$ is faithful for all
$n\geq 1$. \label{Theorem:bikra}
\end{Theorem}
\begin{Remark} Actually,  S~Bigelow \cite{Big2} proved this theorem for 
the representation $\rho_n$ characterized in homological terms and 
D.~Krammer \cite{Kr2} proved the following.
Let $K={\R}  [t^{\pm 1}]$, $q\in {\R}$, and $0<q<1$. Then the 
representation $\rho_n$ defined by (\ref{eq:lakra})  is faithful for all $n\geq 1$.
This result implies Theorem~\ref{Theorem:bikra}:   if a representation over 
${{\Z}}[q^{\pm 1}, t^{\pm 1}]$ becomes faithful  after   assigning a real 
value to $q$, then it is faithful itself. 
\end{Remark} 
M.~G.~Zinno \cite{Zi} established connection between 
Birman-Murakami-Wenzl algebra \cite{BW}, \cite{Mu} and the  
Lawrence-Krammer 
representation. Namely, he proved that the  Lawrence-Krammer 
representation is identical to the irreducible representations of
the Birman-Murakami-Wenzl algebra parametrized by Young diagrams of shapes 
$(n-2)$ and $(1^{n-2}).$ This means that the Young
diagram in the case considered consists of one row (respectively of one 
column) only, with $n-2$ boxes. It follows that Lawrence-Krammer
representation is irreducible.

\section{Generalizations of Braids\label{sec:gen}}
\subsection{Configuration spaces of manifolds\label{subsection:coman}}

The notion of configuration space of Subsection~\ref{subsection:braco}
can be naturally generalized for a configuration space of a 
manifold as follows.
Let $Y$ be a connected topological manifold and let $W$ be a
finite group acting on $Y$. A point $y \in Y$ is called {\it
regular} if its stabilizer $\{w\in W: wy=y \}$ is trivial, i.e.,
consists only of the unit of the group $W$. The set $\tilde Y$ of
all regular points is open. Suppose that it is connected and
nonempty. The subspace $ORB(Y, W)$ of the space of all orbits
$Orb(Y,W)$ consisting of the orbits of all regular points is
called the {\it space of regular orbits}. There is a free action
of $W$ on $\tilde Y$ and the projection $p:\tilde
Y\rightarrow\tilde Y/W=ORB(Y,W)$ defines a covering. Let us
consider the initial segment of the long exact sequence of this
covering: 
\begin{equation}
1 \rightarrow\pi_1(\tilde Y, y_0)\buildrel p_*
\over\rightarrow \pi_1(ORB(Y, W), p(y_0)) \rightarrow W
\rightarrow 1 .\label{eq:pi1}
\end{equation}
The fundamental group $\pi_1(ORB(Y, W),
p(y_0))$ of the space of regular orbits is called the {\it braid
group of the action of $W$ on $Y$} and is denoted by $Br(Y,W)$. The
fundamental group $\pi_1( \tilde Y, y_0)$ is called the {\it pure
braid group of the action of $W$ on $Y$} and is denoted by $P(Y,W)$.
The spaces $\tilde Y$ and $ORB(Y, W)$ are path connected, so the
pair of these groups is defined uniquely up to isomorphism and we
may omit mentioning the base point $y_0$ in the notations.

For any space $Y$ the symmetric group $\Sigma_m$ acts on the
Cartesian power $Y^m$ of the space $Y$ by the formulas 
(\ref{eq:permut})
 We denote
by $F(Y,m)$ the space of $m$-tuples of pairwise different points
in $Y$: $$F(Y,m)=\{(p_1,...,p_m)\in Y^m:p_i\not=p_j \ \text {for} \
i\not=j \}.$$ This is the space of regular points of this action.
In the case when $Y$ is a connected topological manifold $M$
without boundary and $\operatorname{dim} M \geq 2$, the space of
regular orbits $ORB(M^m, \Sigma_m)$ is open, connected and
nonempty. We call $ORB(M^m, \Sigma_m)$ the {\it configuration
space of the manifold} $M$ and denote by $B(M,m)$. The braid group
$Br(M^m, \Sigma_m)$ is called the {\it braid group on $m$ strings
of the manifold $M$} and is denoted by $Br(m,M)$. Analogously, we
call the group $P(M^m, \Sigma_m)$ the {\it pure braid group on $m$
strings of the manifold $M$} and denote it by $P(m,M)$. These
definitions of braid groups were given by R.~Fox and L.~Neuwirth
\cite{FoN}. 

\subsection{Artin-Brieskorn braid groups\label{subsection:babg}}

The braid groups are included in the series of so called
generalized braid groups (this was their name in the work of
E.~Brieskorn of 1971 \cite{Bri1}),
or Artin groups (so they we called by E.~Brieskorn and K.~Saito in 
the paper of 1972 \cite{BS}). They were defined by 
E.~Brieskorn \cite{Bri1}, so we call them Artin--Brieskorn groups.

Let $V$ be a finite dimensional real vector space
($\operatorname{dim} V= n$) with Euclidean structure. Let $W$ be a
finite subgroup of $GL(V)$ generated by reflections. Let
$\mathcal M$ be the set of hyperplanes such that $W$ is generated by
orthogonal reflections with respect to $M \in \mathcal M$. We suppose
that for every $w\in W$ and every hyperplane $M \in \mathcal M$ the
hyperplane $w(M)$ belongs to $\mathcal M$. 

The group $W$ is generated by the reflections
$w_i=w_i(M_i), i \in I$, satisfying only the following relations
$$ (w_i w_j)^{m_{i,j}}=e, \ i, j \in I,$$ 
where the natural
numbers $m_{i,j}=m_{j,i}$ form the {\it Coxeter matrix} of $W$
from which the {\it Coxeter graph} $\Gamma(W)$ of $W$ is
constructed  \cite{Bo}. We use the following notation of
P.~Deligne \cite{Del}: $\operatorname{prod}(m;x,y)$ denotes the product
$xyxy... $ ($m$ factors). The {\it generalized braid group} (or 
{\it Artin--Brieskorn group}) $Br(W)$
{\it of} $W$ \cite{Bri1}, \cite{Del} is defined as the group with 
generators
$\{s_i, i \in I\}$ and relations: 
$$\operatorname{prod}
(m_{i,j};s_i,s_j)= \operatorname{prod} (m_{j,i};s_j,s_i).$$ 
From
this we obtain the presentation of the group $W$ by adding the
relations: 
$$s_i^2=e; \ \  i \in I.$$ 
We see in  
Theorem~\ref{Theorem:bri1} that
this definition of the generalized braid group agrees with our
general definition of a braid group of an action of a group $W$
(Subsection~\ref{subsection:coman}). We denote by $\tau_W$ the canonical homomorphism
from $Br(W)$ to $W$. The classical braids on $k$ strings $Br_k$
are obtained by this construction if $W$ is the symmetric group on
$k+1$ symbols. In this case $m_{i,i+1}=3$, and $m_{i,j}=2$ if $j
\neq i, i+1 $.

Classification of {\it irreducible} (with connected Coxeter graph)
Coxeter groups is well known (see for example Theorem 1 , Chapter
VI, \S 4 of [Bo]). It consists of the three infinite series: $A$,
$C$ (which is also denoted by $B$ because in the corresponding
classification of simple Lie algebras two different series $B$ and
$C$ have this group as their Weyl group) and $D$ as well as the
exceptional groups $ E_6, E_7, E_8, F_4$, $G_2$, $H_3$, $H_4$ and
$I_2(p).$ 

Now let us consider the complexification $V_C$ of the space $V$
and the complexification $M_C$ of $M\in\mathcal M$. Let
$Y_W=V_C-\bigcup_{M\in\mathcal M}M_C$. The group $W$ acts freely 
on $Y_W$. Let $X_W=Y_W/W$
then $Y_W$ is a covering over $X_W$ corresponding to the group
$W$. Let $y_0 \in A_0 $ be a point in some chamber $A_0$ and let
$x_0$ stand for its image in $X_W$. We are in the situation
described in Subsection~\ref{subsection:coman} in the definition 
of the braid group of
the action of the group $W$. This braid group is defined as the
fundamental group of the space of regular orbits of the action of
$W$. In our case $ORB(V_C,W)=X_W.$ So, the generalized braid group
is $\pi_1(X_W,x_0)$. For each $j\in I$, let $\ell_j^\prime$ be the
homotopy class of paths in $Y_W$ starting from $y_0$ and ending in
$w_j(y_0)$ which contains a polygon line with successive vertices:
$y_0, y_0+iy_0, w_j(y_0)+iy_0, w_j(y_0)$. The image $\ell_j$ of
the class $\ell_j^\prime$ in $X_W$ is a loop with base point
$x_0$.

\begin{Theorem} The fundamental group $ \pi_1(X_W,x_0)$
is generated by the elements $\ell_j$ satisfying the following
relations: $$\operatorname{prod} (m_{j,k};\ell_j , \ell_k )=
\operatorname{prod} (m_{k,j}; \ell_k , \ell_j).$$
\label{Theorem:bri1}
\end{Theorem}
This theorem was proved by E.~Brieskorn \cite{Bri2}.

The word problem and the conjugacy problem for Artin-Brieskorn groups
were solved by E.~Brieskorn and K.~Saito \cite{BS} and P.~Deligne
\cite{Del}. 
The biautomatic structure of these groups
was established by R.~Charney \cite{Cha1}. 

In the case when $V$ is complex finite dimensional space and 
$W$ is a
finite subgroup of $GL(V)$ generated by {\it pseudo-reflections}
the corresponding braid groups were studied by M.~Brou\'e, G.~Malle and 
R.~Rouquier \cite{BMR} and also by D.~Bessis and J.~Michel \cite{BM}.  

\subsection{Braid groups of surfaces}
Braid groups of a sphere $Br_n(S^2)$ also have simple geometric 
interpretation
as a group of isotopy classes of braids lying in a layer between two 
concentric spheres. It has the 
presentation with generators $\delta_i$, $i=1, ..., n-1$, and relations: 
\begin{equation}
 \begin{cases} \delta_i \delta_j &=\delta_j \, \delta_i, \ \
\text{if} \ \ |i-j|
>1,\\ 
\delta_i \delta_{i+1} \delta_i &= \delta_{i+1} \delta_i \delta_{i+1}, \\
\delta_1 \delta_2 \dots \delta_{n-2}\delta_{n-1}^2\delta_{n-2} \dots
\delta_2\delta_1 &=1.
\end{cases} \label{eq:sphrelations}
\end{equation}
This presentation was find by O.~Zariski \cite{Za1} in 1936 and then
rediscovered by E.~Fadell and J.~Van Buskirk \cite{FaV} in 1961.

Presentations of braid groups of all closed surfaces were obtained by
G.~P.~Scott \cite{Sc} but look rather complicated.

\subsection{Braid groups in handlebodies}

The subgroup $Br_{1,n+1}$ of the braid group $Br_{n+1}$ consisting
of braids with the fixed first string can be interpreted also as the 
braid group in a solid torus.
Here we study braids in a handlebody of the arbitrary genus $g$.

Let $H_g$ be a handlebody of  genus $g$. The braid group $Br^g_n$
on $n$ strings in $H_g$ was first considered by A.~B.~Sossinsky
\cite{So}. Let $Q_g$ denote a subset of the complex plain
$\Bbb C$, consisting of $g$ different points, 
$Q_g=\{z_1^0, ...,z_g^0\}$, say, $z_i^0=i$. The interior
of the handlebody $H_g$ may be interpreted as the direct product
of the complex plain $\Bbb C$ without $g$ points: $\Bbb C\setminus
Q_g$, and the open interval, for example, $(-1,1)$: 
$$\dot H_g =
(\Bbb C\setminus Q_g) \times (-1, 1).$$ 
The space
$F(\Bbb C\setminus Q_g,n)$ can be interpreted as the complement of
the arrangement of hyperplanes in $\Bbb C^{g+n}$ given by the
formulas: 
$$H_{j,k}: z_j-z_k=0 \ \text{for all} \ j, k;$$ 
$$H^i_j:
z_j=z_i^0 \ \text{for} \ i=1, ..., g; \ j= 1, ...,n.$$
The braids in $Br^g_n$ are
considered as lying between the planes with coordinates $z=0$ and
$z=1$ and connecting the points $((g+1,0),..., (g+n,0))$. So
$Br_n^g$ can be considered as a subgroup of the classical braid
group $Br_{g+n}$ on $g+n$ strings such that the braids from
$Br_n^g$ leave the first $g$ strings unbraided. In this subsection we
denote by $\bar\sigma_j$ the standard generators of the group
$Br_{g+n}$. Let $\tau_k, \ k=1,2,..., g$, be the following braids:
$$\tau_k = \bar\sigma_g\bar\sigma_{g-1}...\bar\sigma_{k+1}
\bar\sigma_{k}^2\bar\sigma_{k+1}^{-1}...\bar\sigma_{g-1}^{-1}
\bar\sigma_g^{-1}.$$ Such a braid is depicted in 
Figure~\ref{fi:tau1}.
\begin{figure}
\begin{picture}(0,130)(-30,0)
\put(-100,100){\line(0,-1){57}} \put(-100,33){\line(0,-1){33}}
\put(-80,100){\line(0,-1){100}} \put(-45,100){\line(0,-1){100}}
\put(-45,85){\line(0,-1){85}} \put(-25,100){\line(-2,-1){15}}
\put(-85,70){\line(-2,-1){10}} \put(-50,88){\line(-2,-1){25}}
\put(-75,25){\line(2,-1){25}} \put(-40,7){\line(2,-1){15}}
\put(-5,100){\line(0,-1){100}} \put(30,100){\line(0,-1){100}}
\put(-105,60){\line(-2,-1){20}} \put(-125,50){\line(2,-1){30}}
\put(-95,35){\line(2,-1){10}}
\put(-100,110){\makebox(0,0)[cc]{$1$}}
\put(-80,110){\makebox(0,0)[cc]{$2$}}
\put(-45,110){\makebox(0,0)[cc]{$g$}}
\put(-25,110){\makebox(0,0)[cc]{$g+1$}}
\put(30,110){\makebox(0,0)[cc]{$g+n$}}
\put(-61,0){\makebox(0,0)[cc]{...}}
\put(13,0){\makebox(0,0)[cc]{...}}
\end{picture}
\caption{}\label{fi:tau1}
\end{figure} 
The elements $\tau_k, \ k=1,2,..., g$, generate a free subgroup
$F_g$ in the braid group $Br_{g+n}$. It
follows for example from the Markov normal form that the elements 
$\tau_k, \ k=1,2,...,g$, together with the standard generators $\bar\sigma_{g+1},
...,\bar\sigma_{g+n-1}$ generate the group $Br^g_n$. So, the braid
group in the handlebody $Br^g_n$ can be considered as a subgroup
of $Br_{g+n}$, generated by two subgroups: $F_g$ and $Br_n$.
Denote by $\sigma_{1}, ..., \sigma_{n-1}$ the standard generators
of $Br_n$ considered as the elements of $Br_n^g$,
$\sigma_{i}=\bar\sigma_{g+i}$, $i=1, \dots, n-1$. So we have the 
presentation of 
$Br_n^g$ with the generators $\tau_k$ and $\sigma_i$ and relations 
\cite{So}, \cite{Ve1}, \cite{Ve4}: 
\begin{equation}\begin{cases}
\sigma_i\sigma_j &=\sigma_j\sigma_i \ \
\text{if} \ \ |i-j| >1,\\ 
\sigma_i \sigma_{i+1} \sigma_i &=
\sigma_{i+1} \sigma_i \sigma_{i+1},\\ 
\tau_k\sigma_i
&=\sigma_i\tau_k \ \ \text{if} \ \ k\geq 1, \ i\geq 2,\\
\tau_k\sigma_1\tau_k\sigma_1&=\sigma_1\tau_k\sigma_1\tau_k, \
k=1,2,..., g,\\
 \tau_k\sigma_1^{-1}\tau_{k+l}\sigma_1 &=
\sigma_1^{-1}\tau_{k+l}\sigma_1\tau_k, \ k=1,2,..., g-1; \
l=1,2,..., g-k. \\
\end{cases} \label{eq:rehan}
\end{equation}
The relation of the fourth type in (\ref{eq:rehan}) is the relation
of the braid group of type $B$ $(C)$.
 The
relations of the fifth type in (\ref{eq:rehan}) describe the interaction
between the generators of the free group and their closest
neighbour $\sigma_1$. Geometrically this is seen in 
Figure~\ref{fi:grehan}. 
\begin{figure}
\begin{picture}(0,270)(0,-50)
\put(-120,200){\line(0,-1){30}} \put(-90,200){\line(1,-1){30}}
\put(-60,200){\line(-1,-1){12}} \put(-30,200){\line(0,-1){120}}
\put(50,200){\line(0,-1){150}} \put(80,200){\line(0,-1){30}}
\put(110,200){\line(1,-1){12}} \put(140,200){\line(-1,-1){42}}
\put(-78,182){\line(-1,-1){24}} \put(128,182){\line(1,-1){12}}
\put(-120,170){\line(1,-1){30}} \put(-60,170){\line(0,-1){60}}
\put(80,170){\line(1,-1){30}} \put(140,170){\line(0,-1){60}}
\put(-108,152){\line(-1,-1){12}} \put(92,152){\line(-1,-1){12}}
\put(-120,140){\line(1,-1){42}} \put(-90,140){\line(-1,-1){12}}
\put(80,140){\line(1,-1){60}} \put(110,140){\line(-1,-1){12}}
\put(-108,122){\line(-1,-1){12}} \put(92,122){\line(-1,-1){12}}
\put(-60,110){\line(-1,-1){30}} \put(-120,110){\line(0,-1){150}}
\put(80,110){\line(0,-1){30}} \put(140,110){\line(-1,-1){12}}
\put(-72,92){\line(1,-1){24}} \put(122,92){\line(-1,-1){24}}
\put(-30,80){\line(-1,-1){42}} \put(-90,80){\line(0,-1){30}}
\put(80,80){\line(1,-1){30}} \put(140,80){\line(0,-1){120}}
\put(-42,62){\line(1,-1){12}} \put(92,62){\line(-1,-1){24}}
\put(-30,50){\line(0,-1){60}} \put(-90,50){\line(1,-1){30}}
\put(50,50){\line(1,-1){30}} \put(110,50){\line(0,-1){60}}
\put(-78,32){\line(-1,-1){12}} \put(62,32){\line(-1,-1){12}}
\put(-90,20){\line(1,-1){60}} \put(-60,20){\line(-1,-1){12}}
\put(50,20){\line(1,-1){42}} \put(80,20){\line(-1,-1){12}}
\put(-78,2){\line(-1,-1){12}} \put(62,2){\line(-1,-1){12}}
\put(-90,-10){\line(0,-1){30}} \put(-30,-10){\line(-1,-1){12}}
\put(50,-10){\line(0,-1){30}} \put(110,-10){\line(-1,-1){30}}
\put(-48,-28){\line(-1,-1){12}} \put(98,-28){\line(1,-1){12}}
\put(10,90){\makebox(0,0)[cc]{$-$}}
\put(10,85){\makebox(0,0)[cc]{$-$}}
\end{picture}
\caption{}\label{fi:grehan}
\end{figure} 
If we introduce the new generators $\theta_k$,  $k=1,2,..., g-1;$
by the formulas:
$$ \theta_k = \sigma_1^{-1}\tau_{k}\sigma_1$$
we obtain the ``positive'' presentation of the group $B_n^g$
with generators of the types $\sigma_i$, $\tau_k$, $\theta_k$
and relations:
\begin{equation}\begin{cases}
\sigma_i\sigma_j &=\sigma_j\sigma_i \ \
\text{if} \ \ |i-j| >1,\\ 
\sigma_i \sigma_{i+1} \sigma_i &=
\sigma_{i+1} \sigma_i \sigma_{i+1},\\ \tau_k\sigma_i
&=\sigma_i\tau_k \ \ \text{if} \ \ k\geq 1, \ i\geq 2,\\
\tau_k\sigma_1\tau_k\sigma_1&=\sigma_1\tau_k\sigma_1\tau_k, \
k=1,2,..., g,\\
 \tau_k\theta_{k+l} &=
\theta_{k+l}\tau_k, \ k=1,2,..., g-1; \
l=1,2,..., g-k. \\
\sigma_1\theta_{k} &=
\tau_k\sigma_1, \ k=1,2,..., g-1. \\
\end{cases} \label{eq:rehanp}
\end{equation}

There is a version of Markov Theorem~\ref{Theorem:mnf}  for
the group $Br^g_n$ \cite{Ve1}, \cite{Ve4}. 
 \subsection{Braids with singularities}

Let $BP_n$ be the subgroup of $\operatorname{Aut} F_n$, generated
by both sets of the automorphisms $\sigma_i$ of (\ref{eq:autf})
and  $\xi_i$ of the following form: 
\begin{equation*}
\begin{cases} x_i &\mapsto x_{i+1}, \\ x_{i+1} &\mapsto
x_i,     \\ x_j &\mapsto x_j, j\not=i,i+1, \end{cases}
\end{equation*}
This is the {\it braid-permutation group}. R.~Fenn,
R.~Rim\'anyi and C.~Rourke proved \cite{FRR1}, \cite{FRR2} that 
this group
is given by the set of generators: $\{ \xi_i, \sigma_i, \ \ i=1,2,
..., n-1 \}$ and relations: 
$$ \begin{cases} \xi_i^2&=1, \\ \xi_i \xi_j
&=\xi_j \xi_i, \ \ \text {if} \ \ |i-j| >1,\\ \xi_i \xi_{i+1}
\xi_i &= \xi_{i+1} \xi_i \xi_{i+1}. \end{cases} $$ 
\vglue 0.01cm
\centerline { The symmetric group relations} 
$$ \begin{cases} \sigma_i \sigma_j &=\sigma_j
\sigma_i, \ \text {if} \  |i-j| >1,
\\ \sigma_i \sigma_{i+1} \sigma_i &= \sigma_{i+1} \sigma_i \sigma_{i+1}.
\end{cases} $$
\vglue 0.01cm
\centerline {The braid group relations } $$ \begin{cases} \sigma_i \xi_j
&=\xi_j \sigma_i, \ \text {if} \  |i-j| >1,
\\ \xi_i \xi_{i+1} \sigma_i &= \sigma_{i+1} \xi_i \xi_{i+1},
\\ \sigma_i \sigma_{i+1} \xi_i &= \xi_{i+1} \sigma_i \sigma_{i+1}.
\end{cases} $$
\vglue0.01cm
\centerline {The mixed relations }
\smallskip

R.~Fenn, R.~Rim\'anyi and C.~Rourke also gave a geometric
interpretation of $BP_n$ as a group of {\it welded braids}. At
first they defined a {\it welded braid diagram} on $n$ strings as
a collection of $n$ monotone arcs starting from $n$ points at a
horizontal line of a plane (the top of the diagram) and going down
to $n$ points at another horizontal line (the bottom of the
diagram). The diagrams may have crossings of two types: 1) the
same as usual braids as for example on the Figure~\ref{fi:sigma} or 
2)~welds as depicted in Figure~\ref{fi:weld}
\begin{figure}
\begin{picture}(0,50)(0,-5) 
\put(0,20){\circle*{5}} \put(-10,40){\line(1,-2){20}}
\put(10,40){\line(-1,-2){20}}
\end{picture}
\caption{}\label{fi:weld}
\end{figure} 

Composition of welded braid diagrams on $n$ strings is defined by
stacking one diagram under the other. The diagram with no
crossings or welds is the identity with respect to composition. So
the set of welded braid diagrams on $n$ strings forms a semigroup
which is denoted by $WD_n$.

R.~Fenn, R.~Rim\'anyi and C.~Rourke defined the 
allowable moves on welded braid diagrams. They 
consist of usual Reidemeister moves (Figure~\ref{fi:R23}) and the 
specific moves depicted of Figures~\ref{fi:R23w}, \ref{fi:Rm},
\ref{fi:Rc}
 \begin{figure}
\begin{picture}(0,110)(0,-15)
\put(-110,60){\circle*{5}} \put(10,60){\circle*{5}}
\put(110,60){\circle*{5}} \put(20,40){\circle*{5}}
\put(100,40){\circle*{5}} \put(-110,20){\circle*{5}}
\put(10,20){\circle*{5}} \put(110,20){\circle*{5}}
\put(-120,80){\line(1,-2){20}} \put(-100,80){\line(-1,-2){20}}
\put(-70,80){\line(0,-1){80}} \put(-50,80){\line(0,-1){80}}
\put(0,80){\line(1,-2){40}} \put(20,80){\line(-1,-2){20}}
\put(40,80){\line(-1,-2){40}} \put(80,80){\line(1,-2){40}}
\put(100,80){\line(1,-2){20}} \put(120,80){\line(-1,-2){40}}
\put(-120,40){\line(1,-2){20}} \put(-100,40){\line(-1,-2){20}}
\put(0,40){\line(1,-2){20}} \put(120,40){\line(-1,-2){20}}
\put(-85,40){\makebox(0,0)[cc]{$\leftrightarrow$}}
\put(60,40){\makebox(0,0)[cc]{$\leftrightarrow$}}
\end{picture}

\caption{}\label{fi:R23w}
\end{figure} 
 \begin{figure}
\begin{picture}(0,110)(0,-15)
\put(-120,60){\circle*{5}} \put(120,60){\circle*{5}}
\put(-110,40){\circle*{5}} \put(-40,40){\circle*{5}}
\put(-30,20){\circle*{5}} \put(30,20){\circle*{5}}
\put(-130,80){\line(1,-2){40}} \put(-120,80){\line(0,-1){80}}
\put(-90,80){\line(-1,-2){27}} \put(-60,80){\line(1,-2){40}}
\put(-30,80){\line(0,-1){80}} \put(-20,80){\line(-1,-2){7}}
\put(20,80){\line(1,-2){40}} \put(30,80){\line(0,-1){13}}
\put(60,80){\line(-1,-2){17}} \put(90,80){\line(1,-2){40}}
\put(120,80){\line(0,-1){53}} \put(130,80){\line(-1,-2){17}}
\put(-130,0){\line(1,2){7}} \put(-60,0){\line(1,2){27}}
\put(20,0){\line(1,2){17}} \put(30,0){\line(0,1){53}}
\put(90,0){\line(1,2){17}} \put(120,0){\line(0,1){13}}
\put(-75,40){\makebox(0,0)[cc]{$\leftrightarrow$}}
\put(75,40){\makebox(0,0)[cc]{$\leftrightarrow$}}
\end{picture}
\caption{}\label{fi:Rm}
\end{figure} 
 \begin{figure}
\begin{picture}(0,110)(0,-15)
\put(40,60){\circle*{5}} \put(-90,20){\circle*{5}}
\put(-100,80){\line(0,-1){40}} \put(-80,80){\line(0,-1){40}}
\put(-50,80){\line(1,-2){20}} \put(-30,80){\line(-1,-2){8}}
\put(30,80){\line(1,-2){20}} \put(50,80){\line(-1,-2){20}}
\put(80,80){\line(0,-1){40}} \put(100,80){\line(0,-1){40}}
\put(-100,40){\line(1,-2){20}} \put(-80,40){\line(-1,-2){20}}
\put(-50,40){\line(1,2){8}} \put(-50,40){\line(0,-1){40}}
\put(-30,40){\line(0,-1){40}} \put(30,40){\line(0,-1){40}}
\put(50,40){\line(0,-1){40}} \put(80,40){\line(1,-2){20}}
\put(100,40){\line(-1,-2){8}} \put(80,0){\line(1,2){8}}
\put(0,40){\makebox(0,0)[cc]{$\leftrightarrow$}}
\put(-110,40){\makebox(0,0)[cc]{$. . .$}}
\put(-65,40){\makebox(0,0)[cc]{$. . .$}}
\put(-20,40){\makebox(0,0)[cc]{$. . .$}}
\put(20,40){\makebox(0,0)[cc]{$. . .$}}
\put(65,40){\makebox(0,0)[cc]{$. . .$}}
\put(110,40){\makebox(0,0)[cc]{$. . .$}}
\end{picture}
\caption{}\label{fi:Rc}
\end{figure} 
The automorphisms of $F_n$ which lie in $BP_n$ can be
characterized as follows. Let $\pi\in\Sigma_n$ be a
permutation and $w_i, \ \ i=1,2...,n$, be words in $F_n$. Then the
mapping $$x_i\mapsto w_i^{-1}x_{\pi(i)}w_i$$ determines an
injective endomorphism of $F_n$. If it is also surjective, we call
it an automorphism of {\it permutation-conjugacy type}. The
automorphisms of this type comprise a subgroup of
$\operatorname{Aut} F_n$ which is precisely $BP_n$. 

The {\it Baez--Birman monoid} $SB_n$  or {\it singular braid monoid} 
\cite{Bae}, \cite{Bir2} is defined as the
monoid with generators $g_i,g_i^{-1},a_i$, $i=1,\dots,n-1,$ and
relations 
\begin{eqnarray*}
&g_ig_j=g_jg_i, \ \text {if} \ \ |i-j| >1,\\
&a_ia_j=a_ja_i, \ \text {if} \ \ |i-j| >1,\\ &a_ig_j=g_ja_i, \
\text {if} \ \ |i-j| \not=1,\\ &g_i g_{i+1} g_i = g_{i+1} g_i
g_{i+1},\\ &g_i g_{i+1} a_i = a_{i+1} g_i g_{i+1},\\ &g_{i+1} g_i
a_{i+1} = a_i g_{i+1} g_i,\\ &g_ig_i^{-1}=g_i^{-1}g_i =1.
\end{eqnarray*}

In pictures $g_i$ corresponds to canonical generator of the braid
group and $a_i$ represents an intersection
of the $i$th and $(i+1)$st strand as in 
Figure~\ref{fi:singen}. More detailed geometric interpretation of the
Baez--Birman monoid can be found in the article of J.~Birman
\cite{Bir2}.
 \begin{figure}
\begin{picture}(0,130)(0,-10) 
\put(0,50){\circle*{5}} \put(-100,100){\line(0,-1){100}}
\put(-50,100){\line(0,-1){100}} \put(-25,100){\line(1,-2){50}}
\put(25,100){\line(-1,-2){50}} \put(50,100){\line(0,-1){100}}
\put(100,100){\line(0,-1){100}}
\put(-100,110){\makebox(0,0)[cc]{$1$}}
\put(-50,110){\makebox(0,0)[cc]{$i-1$}}
\put(-25,110){\makebox(0,0)[cc]{$i$}}
\put(25,110){\makebox(0,0)[cc]{$i+1$}}
\put(50,110){\makebox(0,0)[cc]{$i+2$}}
\put(100,110){\makebox(0,0)[cc]{$n$}}
\put(-75,50){\makebox(0,0)[cc]{.\quad.\quad.}}
\put(75,50){\makebox(0,0)[cc]{.\quad.\quad.}}
\end{picture}
\caption{}\label{fi:singen}
\end{figure} 
 R. Fenn, E. Keyman and C. Rourke
proved \cite{FKR} that the Baez--Birman monoid embeds in a group $SG_n$
which they called the {\it singular braid group}: $$SB_n\to
SG_n.$$ So, in $SG_n$ the elements $a_i$ become invertible and all
relations of $SB_n$ remain true.

\begin{figure}
\epsfbox{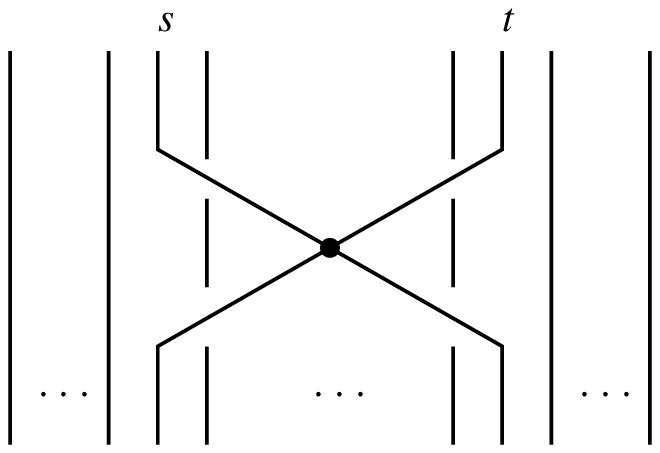}
\caption{}
\label{fi:sbige3}
\end{figure}

The analogue of the Birman-Ko-Lee presentation for the singular braid monoid 
was obtained  in \cite{V10}. Namely, it was proved that the  
monoid $SB_n$ has a presentation with
generators $a_{ts}$, $a_{ts}^{-1}$
for $1\leq s<t\leq n$ and  $b_{qp}$ for
$1\leq p<q\leq n$
and relations
\begin{equation} \begin{cases}
a_{ts}a_{rq}&=a_{rq}a_{ts} \ \ {\rm for} \ \ (t-r)(t-q)(s-r)(s-q)>0,\\
a_{ts}a_{sr} &=a_{tr}a_{ts}=a_{sr}a_{tr}  \ \ {\rm for} \ \
1\leq r<s<t\leq n , \\
a_{ts}a_{ts}^{-1} &=a_{ts}^{-1}a_{ts} =1 \ \ {\rm for} \ \ 1\leq s<t\leq n,\\
a_{ts}b_{rq}&=b_{rq}a_{ts} \ \ {\rm for} \ \ (t-r)(t-q)(s-r)(s-q)>0,\\
a_{ts}b_{ts}&=b_{ts}a_{ts}  \ \ {\rm for} \ \
1\leq s<t\leq n , \\
a_{ts}b_{sr} &=b_{tr}a_{ts}  \ \ {\rm for} \ \
1\leq r<s<t\leq n , \\
a_{sr}b_{tr} &=b_{ts}a_{sr}  \ \ {\rm for} \ \
1\leq r<s<t\leq n , \\
a_{tr}b_{ts}&=b_{sr}a_{tr}  \ \ {\rm for} \ \
1\leq r<s<t\leq n, \\
b_{ts}b_{rq}&=b_{rq}b_{ts} \ \ {\rm for} \ \ (t-r)(t-q)(s-r)(s-q)>0.
\end{cases}\label{eq:srebkl}
\end{equation}
The elements  $a_{ts}$ are defined the same way as in \ref{eq:ats}
and  the elements $b_{qp}$ for 
$1\leq p<q\leq n$  are
defined by
\begin{equation} 
b_{qp}=(\sigma_{q-1}\sigma_{q-2}\cdots\sigma_{p+1}) x_p
(\sigma^{-1}_{p+1}\cdots\sigma^{-1}_{q-2}\sigma^{-1}_{q-1})  \ \
{\rm for} \ \ 1\leq p<q\leq n.
\label{eq:defab}
\end{equation}
Geometrically the generators  $b_{s,t}$ are depicted
in Figure~\ref{fi:sbige3}.

\section{Homological Properties\label{sec:hom}}
\subsection{Configuration spaces and $K(\pi,1)$-spaces}

Let $(q_i)_{i \in \Bbb N}$ be a fixed sequence of distinct points
in the manifold $M$ and put $Q_m =\{q_1, ..., q_m\}$. We use
$$Q_{m,l}= (q_{l+1},...,q_{l+m}) \in F(M\setminus Q_l,m)$$ as the
standard base point of the space $F(M\setminus Q_l,m)$. If $k<m$
we define the projection $$\operatorname{proj} : F(M\setminus
Q_l,m) \rightarrow F(M\setminus Q_l,k)$$ by the formula:
$\operatorname{proj}(p_1,...,p_m) = (p_1,...,p_k).$ The following
theorems were proved by E.~Fadell and L.~Neuwirth \cite{FaN}.
\begin{Theorem} The triple $\operatorname{proj}:
F(M\setminus Q_l,m)\rightarrow F(M\setminus Q_l,k)$ is a locally
trivial fiber bundle with fibre $\operatorname{proj}^{-1}Q_{k,l}$
homeomorphic to $F(M\setminus Q_{k+l},m-k)$.
\label{Theorem:fan1}
\end{Theorem}
Consideration of the sequence of fibrations $$F(M\setminus
Q_{m-1},1) \rightarrow F(M\setminus Q_{m-2},2) \rightarrow
M\setminus Q_{m-2},$$ $$F(M\setminus Q_{m-2},2) \rightarrow
F(M\setminus Q_{m-3},3) \rightarrow M\setminus Q_{m-3},$$ $$... \
,$$ $$F(M\setminus Q_{1},m-1) \rightarrow F(M,m) \rightarrow M$$
leads to the following theorem. 
\begin{Theorem} For any
manifold $M$ $$\pi_i (F(M\setminus Q_1,m-1)) =
\oplus_{k=1}^{m-1}\pi_i (M\setminus Q_k)$$ for $i \geq 2$. If
$\operatorname{proj} : F(M,m) \rightarrow M$ admits a section then
$$\operatorname{proj}_i\pi_i(F(M,m))=\oplus_{k=0}^{m-1}\pi_i(M\setminus
Q_k), \ i \geq 2.$$
\label{Theorem:fan2}
\end{Theorem}
\begin{Corollary} If $M$ is the Euclidean $r$-space, then
$$\pi_i(F(M,m))= \oplus_{k=0}^{m-1}\pi_i(\underbrace{S^{r-1}\vee
...\vee S^{r-1}}_k), \ i\geq 2.$$
\end{Corollary}
\begin{Corollary} If $M$ is the Euclidean $2$-space, then
$F(\Bbb R^2, m)$ is the $K(P_m, 1)$-space and $B(\Bbb R^2, m)$ is
the $K(Br_m, 1)$-space.
\label{Corollary:fan}
\end{Corollary}

Let $X_W$ be the space defined in Subsection \ref{subsection:babg}.
\begin{Theorem} The universal covering of $X_W$ is contractible, 
and so $X_W$ is a $K( \pi;1)$-space.
\label{Teorem:kapi1}
\end{Theorem}
This theorem for the groups of types $C_n, \ G_2$ and $I_2(p)$,
was proved by E.~Brieskorn \cite{Bri1} in much the same way as
E.~Fadell and L.~Neuwirth \cite{FaN} proved Theorems~\ref{Theorem:fan1},
\ref{Theorem:fan2} and
Corollary~\ref{Corollary:fan}.  For the groups
of types $D_n$ and $F_4$ E.~Brieskorn used this method with minor
modifications. In general case  Theorem~\ref{Teorem:kapi1} was
proved by P.~Deligne \cite{Del}.

It follows from Theorem~\ref{Theorem:fan2}  that $F(\Bbb C\setminus Q_g,n)$ and
$B(\Bbb C\setminus Q_g)$ are $K(\pi,1)$-spaces, \penalty-10000
$\pi_1B(\Bbb C\setminus Q_g)=Br_n^g$, so, $B(\Bbb C\setminus Q_g)$
can be considered as the classifying space of $Br_n^g$. 

\subsection{Cohomology of pure braid groups}

Cohomology of pure braid groups were first calculated by
V.~I.~Arnold \cite{Arn1}. The map $$ \phi: S^{n-1} \rightarrow F(\Bbb R^n, 2),$$
described by the formula $\phi(x)= (x,-x)$, is a
$\Sigma_2$-equivariant homotopy equivalence. Denote by $A$ the
generator of $H^{n-1}(F(\Bbb R^n, 2), \Bbb Z)$ that is mapped by
$\phi^*$ to the standard generator of $H^{n-1}(S^{n-1},\Bbb Z)$.
For $i$ and $j$, such that $1\leq i,j\leq m, \ i\not= j$, specify
$\pi_{i,j}:F(\Bbb R^n, m)\rightarrow F(\Bbb R^n, 2)$ by the
formula $\pi_{i,j}(p_1, ..., p_m)= (p_i,p_j)$. Put
$$A_{i,j}=\pi_{i,j}^*(A) \in H^{n-1}(F(\Bbb R^n, m), \Bbb Z).$$ It
follows that $A_{i,j}= (-1)^n A_{j,i}$ and $A_{i,j}^2= 0$. For
$w\in\Sigma_m$ there is an action $w(A_{i,j})=
A_{w^{-1}(i),w^{-1}(j)}$, since $\pi_{i,j} w =
\pi_{w^{-1}(i),w^{-1}(j)}.$ Note also that under restriction to
$$F(\Bbb R^n \setminus Q_k, m-k)\cong\pi^{-1}(Q_k) \subset F(\Bbb
R^n, m),$$ the classes $A_{i,j}$ with $1 \leq i,j \leq k$ go to
zero since in this case the map $\pi_{i,j}$ is constant on
$\pi^{-1}(Q_k)$.
\begin{Theorem} The cohomology group
$H^*(F(\Bbb R^n \setminus Q_k, m-k), \Bbb Z)$ is the free Abelian
group with generators $$A_{i_1,j_1}A_{i_2,j_2} ... A_{i_s,j_s},$$
where $k<j_1<j_2<...<j_s \leq m$ and $i_r<j_r$ for $r=1,\ \dots,
s.$
\end{Theorem}
The multiplicative structure and the $\Sigma_m$-algebra structure
of $H^*(F(\Bbb R^n, m), \Bbb Z)$ are given by the following
theorem which is proved using the $\Sigma_3$-action on $H^*(F(\Bbb
R^n, 3), \Bbb Z)$. 
\begin{Theorem} The cohomology ring
$H^*(F(\Bbb R^n, m), \Bbb Z)$ is multiplicatively generated by the
square-zero elements $$A_{i,j} \in H^{n-1}(F(\Bbb R^n, m), \Bbb
Z), \ 1 \leq i<j\leq m ,$$ subject only to the relations
\begin{equation}
A_{i,k}A_{j,k}= A_{i,j} A_{j,k} - A_{i,j} A_{i,k} \ \text {for}
\ i<j<k. \label{eq:reqp}
\end{equation}
 The Poincar\'e series for $F(\Bbb R^n, m)$ is
the product \ $\prod_{j=1}^{m-1}(1+jt^{n-1}).$
\end{Theorem}
\begin{Remark} In the case of $\Bbb R^2 = \Bbb C$ the cohomology
classes $A_{j,k}$ can be interpreted as the classes of cohomology
of differential forms $$\omega_{j,k}={1\over 2\pi i}
{dz_j-dz_k\over z_j-z_k}.$$
\end{Remark}

E.~Brieskorn calculated the cohomology of pure generalized braid
groups \cite{Bri1} using ideas of V.~I.~Arnold for the classical case.
Let $\mathcal V$ be a finite-dimensional complex vector space and
$H_j\in \mathcal V, \ j\in I$ be the finite family of complex affine
hyperplanes given by linear forms $l_j$. E.~Brieskorn proved the
following fact. 
\begin{Theorem} The cohomology classes,
corresponding to the holomorphic differential forms
$$\omega_j={1\over 2\pi i} {dl_j\over l_j},$$ generate the
cohomology ring $H^*(\mathcal V\setminus \cup_{j\in I}H_j,\Bbb Z).$
Moreover, this ring is isomorphic to the $\Bbb Z$-subalgebra
generated by the forms $\omega_j$ in the algebra of meromorphic
forms on $\mathcal V.$
\end{Theorem}
The cohomologies of pure generalized braid groups are described as
follows. 
\begin{Theorem} (i) The cohomology group
$H^k(P(W), \Bbb Z)$ of the pure braid group $P(W)$ with integer
coefficients is a free Abelian group, and its rank is equal to the
number of elements $w\in W$ of length $l(w)=k$, where $l$ is the
length considered with respect to the system of generators
consisting of all reflections of $W$.

(ii) The Poincar\'e series for $H^*(P(W), \Bbb Z)$ is the product
$\prod_{j=1}^{n}(1+m_jt),$ where $m_j$ are the exponents of the
group $W$.

(iii) The multiplicative structure of $H^*(P(W), \Bbb Z)$
coincides with the structure of algebra, generated by 1-forms
described in the previous theorem.
\end{Theorem}
\subsection{Homology of braid groups}
To study the cohomologies of the classical braid groups $H^*(Br_n,
\Bbb Z)$, V.~I.~Arnold \cite{Arn2} interpreted the space
$K(Br_n,1)\cong B(\Bbb R^2, n)$ as the space of monic complex
polynomials of degree $n$ without multiple roots
$$P_n(t)=t^n+z_1t^{n-1}+...+z_{n-1}t+z_n.  $$
Using this idea he proved 
Theorems of finiteness, of recurrence and of stabilization.
Homology with coefficients in $\Z/2$ were calculated by D.~B.~Fuks
in the following theorems \cite{Fuks1}.
\begin{Theorem} The
homology of the braid group on the infinite number of strings with
coefficients in $\Bbb Z/2$ as a Hopf algebra is isomorphic to the
polynomial algebra on infinitely many generators $a_i, i=1,2,...;$
$\operatorname{deg} \, a_i =2^i-1$: $$H_*(Br_\infty ,\Bbb Z/2)
\cong \Bbb Z/2[a_1, a_2, ..., a_i, ...]$$ with the coproduct given
by the formula: $$\Delta(a_i)=1\otimes a_i + a_i\otimes 1.$$
\end{Theorem}
\begin{Theorem} The canonical inclusion $ Br_n \rightarrow
Br_\infty$ induces a monomorphism in homology with coefficients in
$\Bbb Z/2$. Its image is the subcoalgebra of the polynomial
algebra $\Bbb Z/2[a_1, a_2, ..., a_i, ...]$ with $\Bbb Z/2$-basis
consisting of monomials $$ a_{1}^{k_1} ... a_{l}^{k_l} \
\text{such that} \ \sum_{i} k_i2^i \leq n.$$
\end{Theorem}
\begin{Theorem} The canonical homomorphism $Br_n\rightarrow
BO_n, \ 1 \leq n \leq \infty$ induces a monomorphism (of Hopf
algebras if $n= \infty$) $$H_*(Br_n,\Bbb Z/2)  \rightarrow
H_*(BO_n,\Bbb Z/2). $$
\end{Theorem}
F.~R.~Cohen calculated the homology of braid groups with
coefficients  $\Bbb Z/p, \ p>2$ also as modules over the 
Steenrod algebra \cite{CF1}, \cite{CF2}, \cite{CF3}.

Later V.~V.~Goryunov  \cite{Go1}, \cite{Go2} applied the methods of Fuks
 and expressed the cohomologies of the generalized
braid groups of types $C$ and $D$ in terms of the cohomologies of
the classical braid groups. 

\section{Connections with the Other Domains\label{sec:con}}
\subsection{Markov Theorem}
\begin{figure}
\epsfbox{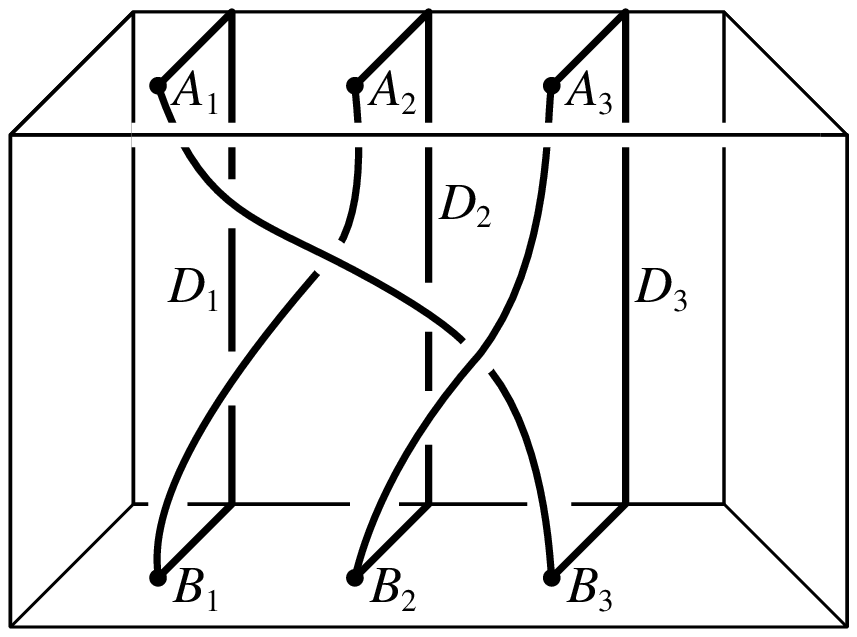}
\caption{} \label{fi:closure}
\end{figure}

Suppose a braid depicted in Figure~\ref{fi:braid} is placed in a cube. On the
boundary of the cube join the point $A_i$ to the point $B_i$ by a mutually 
disjoint simple arc $D_i$. Since our initial braid does not intersect the 
boundary of the cube 
except at the points $A_1, \dots , A_n$ and $B_1, \dots , B_n$ we obtain 
a link
(or, in particular, a knot) i.e. a system of simple closed curves in $\R^3$.  
A link
obtained in such a manner is called the {\it closure} of the braid, see 
Figure~\ref{fi:closure}.

\begin{Theorem} {\rm(J.~W.~Alexander)}
Any link can be represented by a closed braid.
\label{Theorem:Alex}
\end{Theorem}
The next step is to understand equivalence classes of braids which correspond to
links. The following Markov Theorem gives an  answer to this question. At first 
we define two types of Markov moves for braids.

\par{\it Type 1 Markov move} replaces a braid $\beta$ on $n$ strings by its 
conjugate $\gamma \beta \gamma^{-1}$.

\par{\it Type 2 Markov move} replaces a braid $\beta$ on $n$ strings by 
the braid $j_n(\beta) \sigma_n$ on $n+1$ strings or by 
$j_n(\beta) \sigma_n^{-1}$  where  $j_n$ is the canonical inclusion of the group
$Br_n$ into the group $Br_{n+1}$ (see Subsection~\ref{subsec:ini}) 
$$j_n:Br_{n} \to Br_{n+1}.$$

\begin{Theorem} {\rm(A.~A.~Markov)}
Suppose that $\beta$ and $\beta^\prime$ are two braids (not necessary 
with the same number of strings). Then, the closures  of $\beta$ and $\beta^\prime$ 
represent the same link if and only if $\beta$ can be transformed into 
$\beta^\prime$ by means of a finite number of type~1 and type~2 Markov moves.
Namely there exists the following sequence,
$$ \beta = \beta_0\to\beta_1\to\dots\to\beta_m=\beta^\prime,$$ 
such that, for $i=0, 1 \dots, m-1,$ $\beta_{i+1}$ is obtained from  $\beta_i$
by the application of a type~1 or 2 Markov moves or their inverses.
\label{Theorem:Mar2}
\end{Theorem}
In  other words, if we consider the disjoint union of all braid groups 
$$\coprod_{n=1}^{n=\infty} Br_n,$$ 
then the Markov moves of types~1 and 2
define the equivalence relation on this set $\sim$ such that the 
quotient set  
$$\coprod_{n=1}^{n=\infty} Br_n/\sim$$
is in one-to-one correspondence with isotopy classes of links.

There exist a lot of proofs of Markov Theorem, see for example the work
of P.~Traczyk \cite{Tr1}.
\subsection{Homotopy groups of spheres and Makanin braids}
Consider the coordinate projections for the spaces $F(M,m)$ where M is a 
manifold (see Subsection~\ref{subsection:coman})
$$d_i: F(M,n+1)\to F(M,n), \ \ i= 0, \dots, n,$$
defined by the formula
$$d_i(p_1, \dots, p_{i+1}, \dots, p_{n+1})= 
(p_1, \dots, \hat{p}_{i+1}, \dots, p_{n+1}).$$
By taking  the fundamental group the maps $d_i$ induces  group homomorphisms 
$${d_i}_*: P_{m+1}(M)\to P_{m}(M), \ \ i= 0, \dots, n.$$
A braid $\beta \in Br_{n+1}$ is called {\it Makanin}  ({\it smooth} 
in the terminology of 
D.~L.~Johnson \cite{Joh1}, {\it Brunnian} in the terminology of 
J.~A.~Berrick, F.~R.~Cohen, Y.~L.~Wong and J.~Wu \cite{BCWW})  if 
$d_i(\beta)= 1$ for all $0\leq i \leq n$. We call them Makanin, because
up to our knowledge it was G.~S.~Makanin who first mentioned them
\cite{Kou}, page 78, question 6.23.
In other words the group of
Makanin braids $Mak_{n+1}(M)$ is given by the formula
$$Mak_{n+1}(M) = \cap_{i=0}^n 
\operatorname{Ker}({d_i}_*: P_{m+1}(M)\to P_{m}(M)).$$

The canonical embedding of the open disc $D^2$ into the sphere $S^2$
$$f: D^2\to S^2$$
induces a group homomorphism
$$f_*: Mak_{n}(D^2) \to Mak_n(S^2 $$
where $Mak_{n}(D^2) $ is the Makanin subgroup $Mak_{n} $ of the 
classical braid group $Br_n$. The group $Mak_{n} $ is free \cite{Gur1, Joh1}.
The following theorem is proved in \cite{BCWW}.  
\begin{Theorem} The is an exact sequence of groups
$$
1\to Mak_{n+1}(S^2) \to Mak_{n}(D^2) \to Mak_{n}(S^2) \to 
\pi_{n-1}(S^2)\to 1
$$ 
for  $n\geq 5$.
\end{Theorem}
Here as usual $\pi_k(S^2)$ denote the $k$-th homotopy group of the sphere $S^2$.

For instance, $Mak_5(S^2) $ modulo $Mak_5 $ is $\pi_4(S^2)=\Z/2$.
The other homotopy groups of $S^2$ are as follows

$$\pi_5(S^2)=\Z/2, \ \pi_6(S^2)=\Z/12, \ \pi_7(S^2)=\Z/2, \ 
\pi_8(S^2)=\Z/2, \dots .$$
Thus, up to certain range,  $Mak_n(S^2) $ modulo $Mak_n $ are
known by nontrivial calculation of $\pi_*(S^2)$.

\end{document}